\documentclass[12pt,letterpaper]{article}

\usepackage{amsfonts}
\usepackage{amsmath}\numberwithin{equation}{section}
\usepackage{amsthm}
\usepackage{amssymb}
\usepackage{bm}
\usepackage[tmargin=1.0in, lmargin=0.75in, rmargin=0.75in, bmargin=1.0in]{geometry}
\usepackage{hyperref}
\usepackage{mathrsfs}
\usepackage{mathtools}
\usepackage[T1]{fontenc}
\usepackage{xparse}
\usepackage{graphicx}
\usepackage[dvipsnames]{xcolor}
\usepackage{multirow}
\usepackage[all]{xy}

\usepackage{tikz-cd}

\usepackage{pst-node}

\newcommand{\R}{\mathbb{R}}

\let\oldd\d
\renewcommand{\d}{%
  \relax\ifmmode\mathrm{d}\else\oldd\fi
}
\renewcommand{\dh}{%
  \relax\ifmmode\mathrm{d}_{1,0}\else\olddh\fi
}

\NewDocumentCommand{\Cinf}{G{}}{\IfNoValueTF{#1}{C^{\infty}}{C^{\infty}_{#1}}}

\renewcommand{\hat}{\widehat}
\renewcommand{\H}[1]{%
  \ifmmode\mathscr{H}^{#1}\else\oldH{#1}\fi
}
\let\oldH\H

\newtheorem{theorem}{Theorem}[section]

\newtheorem{definition}[theorem]{Definition}
\newtheorem{example}[theorem]{Example}

\newtheorem{lemma}[theorem]{Lemma}

\newtheorem{remark}[theorem]{Remark}

\newtheorem{proposition}[theorem]{Proposition}

\numberwithin{equation}{section}

\title{Jacobi fields and Jacobi sections for semisprays on Lie algebroids: transverse variations and equivalences}

\author{Misael Avenda\~{n}o-Camacho, Jhonny Kama-Mamani}

\begin{document}
\maketitle

\begin{abstract}

We study the relationship between Jacobi sections and Jacobi fields associated with semisprays on Lie algebroids. First, we characterize compatible transverse variations in terms of Lie algebroid morphisms. This characterization leads to a variational description of Jacobi sections: a section along an integral curve of a semispray is a Jacobi section if and only if it arises as the infinitesimal variation of a compatible transverse variation of a dynamical variation satisfying a suitable condition on its image in the base. We also investigate the relationship between Jacobi sections and Jacobi fields through the anchor map. We prove that every Jacobi section projects to a Jacobi field, and establish a converse when the corresponding base variation is contained in a single leaf of the characteristic foliation.
%A key ingredient is the construction of a global homotopy operator, which allows us
%to prove the vanishing of the cohomology associated with the vertical component of
%the Lie algebroid differential and to obtain a structural description of a class
%of symplectic sections. As a consequence, we derive an explicit family of Poisson
%structures with the desired dynamical property.

\end{abstract}

%\keywords{Lie algebroids, semisprays, Jacobi Fields, Jacobi sections,  Lie algebroid morphisms, transverse variations}

\section{Introduction}\label{sec1}
Jacobi fields arise naturally from the study of variations through geodesics on Riemannian manifolds or on manifolds with affine connections. They are the variational vector fields associated with such variations and describe the infinitesimal behavior of geodesics under perturbations of their initial conditions. Equivalently, Jacobi fields can be characterized as vector fields $J$ along any geodesic curve $\gamma$ satisfying the Jacobi equation:
 
$$ D_t^2J + R(J,\dot{\gamma})\dot{\gamma}=0, $$
where $D_t$ is the covariant derivative along $\gamma$ and $R$ denotes the curvature tensor. Jacobi fields play an important role in several fundamental results in Riemannian geometry, including the Rauch comparison theorem, the Myers theorem, and the Cartan–Hadamard theorem.

The Levi-Civita connection on a Riemannian manifold $(M,g)$ induces a spray $\mathcal{G}\in \mathfrak{X}(TM)$ which is called the geodesic spray. Then, any geodesic curve on $M$ is the projection of an integral curve of $\mathcal{G}$ by the canonical tangent bundle projection.  Indeed, there exists a correspondence between variations of geodesics on $M$ and variations by integral curves on  $\mathcal{G}$ and their corresponding infinitesimal variations. Hence, one can extend the study of variations of integral curves and infinitesimal variations of any spray or semispray on the tangent bundle. The semisprays on the tangent bundle are very important objects in mechanics since their integral curves describe the dynamics on the manifold governed by second-order differential equations. In this setting, Buc\u{a}taru \cite{Bucataru-1999} derived the concept of Jacobi fields for a spray on the tangent bundle by proving that the infinitesimal variation of a variation of integral curves satisfies a kind of Jacobi equation which involves the curvature of the nonlinear connection induced by the spray. Later, Cari\~{n}ena et al. \cite{carinena2015jacobi} define the notion of Jacobi field for semisprays on the tangent bundle. Their definition of Jacobi fields is based on the infinitesimal variations of variations of integral curves on the base (i.e., projection on the base manifold of integral curves of the semisprays). They also derive a Jacobi-type equation for Jacobi fields of semisprays in terms of the covariant derivative along the semispray and the generalized Jacobi endomorphism defined by the nonlinear connections induced by the vertical endomorphism on the tangent bundle. Anahory Sim\~{o}es et al. introduce the notion of Jacobi fields and the corresponding Jacobi equation for nonholonomic mechanical systems using nonholonomic connections in \cite{anahory2022jacobi}. 

The notion of Jacobi field has also been extended to Lie algebroids in different contexts, but mostly using a variational approach. For example, Boucetta introduced the notion of \textit{Jacobi sections} along geodesic curves on a Riemannian Lie algebroid, as the sections along a given geodesic curve satisfying a Jacobi-type equation defined by the covariant derivative and the curvature associated with the Levi-Civita connection of the metric (\cite{Boucetta2011}).  Furthermore, Boucetta also proved that his definition of Jacobi sections admits a variational approach (see \cite{Boucetta2011}).  In \cite{Jozwikowski-2013},  J\'o\'zwikowski introduces \textit{Jacobi vector fields} for Lagrangian systems on skew-symmetric algebroids using the Legendre mapping and the Tulczyjew differential. Despite the name, they are sections of the algebroid along a base curve. For semisprays on Lie algebroids, Cari\~{n}ena, Gheorghiu and Mart\'inez \cite{carinena2015jacobi} define \textit{Jacobi sections} along an integral curve of the semispray by means of a special kind of Lie algebroid morphisms and derive for them a Jacobi equation of the form $\nabla\nabla W + \Phi(W)=0$ which is expressed in terms of the dynamical covariant derivative $\nabla$ and the Jacobi endomorphism $\Phi$ of the semispray. Their constructions, including a variational characterization, are carried out on the prolongation. Whether this variational characterization can be reformulated in terms of variations of integral curves on the Lie algebroid itself, as in Boucetta's setting, is the question we address in this paper.

In \cite{Avendano-Kama-2026}, we introduced the notion of Jacobi fields for semisprays on Lie algebroids: we consider a variation by integral curves of a given integral curve of the semispray, then we project the variation to the base manifold and define a Jacobi field as the variational vector field of the projected variation. Our definition recovers the classical definition of Jacobi fields in the case of the Riemannian spray of a Riemannian manifold. Actually, when we restrict to the case of sprays on a general Lie algebroid, we proved that they are in a one-to-one correspondence with the linear connections on Lie algebroids: given a spray, there exists a unique torsionless linear connection such that the integral curves of the spray are precisely the geodesic curves of the linear connection. Then, we can extend Boucetta's definition of a Jacobi section along a geodesic (integral) curve to the case of a linear connection not necessarily given by a Riemannian metric on the Lie algebroid. At this point, we can associate to a given spray on a Lie algebroid its Jacobi sections and its Jacobi fields. We also show that every Jacobi section projects via the anchor map to a Jacobi field and, conversely, every Jacobi field tangent to a leaf of the characteristic foliation is the image of a Jacobi section.

%%%%%%%%%%%%%%%%%%%%%%%%%%%
%%% Nueva versión
%%%%%%%%%%%%%%%%%%%%%%%%%%

The goal of this paper is to relate the Jacobi sections defined by Cari\~{n}ena, Gheorghiu and Mart\'inez in \cite{carinena2015jacobi} with the Jacobi fields introduced in \cite{Avendano-Kama-2026}. Since the constructions of \cite{carinena2015jacobi} are carried out on the prolongation of the Lie algebroid, while the notion of Jacobi fields is defined in terms of variations of integral curves of the semispray on the Lie algebroid itself, our first task is to translate between the two languages. 

First, in Theorem \ref{theo-morph-var} we show that a vector bundle morphism $\Phi=\zeta\,dt+\kappa\,ds$ from $T\tilde I\times TI$ into $A$ is a Lie algebroid morphism if and only if $\zeta$ is a variation by admissible curves, with base contained in a leaf of the characteristic foliation, and $\kappa$ is a compatible transverse variation of $\zeta$. Then, we prove our first main theorem (Theorem~\ref{theo-jacobi-transv}): a section is a Jacobi section along an integral curve of a semispray if and only if it is the infinitesimal variation of a compatible transverse variation of a variation by integral curves of the semispray. This variational characterization is analogous to the one given by Boucetta in \cite{Boucetta2011} for geodesics of Riemannian Lie algebroids.  A technical, but relevant, ingredient of the proof is that the lift of a vector field to the prolongation of the tangent Lie algebroid is a Lie algebroid morphism (Proposition~\ref{prop-CX}).

In our formulation, the compatibility condition is the only ingredient whose definition still relies on the prolongation; we leave open the question of describing it intrinsically on the Lie algebroid.

 The second step is to compare Jacobi sections with Jacobi fields. Our second main result states that the anchor map relates both notions: every Jacobi section projects to a Jacobi field (Theorem \ref{sect-field-jac-theo}), and, conversely, every Jacobi field arising from a variation whose projection to the base lies in a single leaf of the characteristic foliation is the image of a Jacobi section (Theorem \ref{sect-field-jac-theo2}). In \cite{Avendano-Kama-2026} this relation was proved for sprays, where Jacobi sections can be defined through the torsionless linear connections associated with the spray, and it was stated in \cite[Remark 6]{Avendano-Kama-2026} that the same relation holds for semisprays. Here we establish it for arbitrary semisprays, using the Jacobi sections of \cite{carinena2015jacobi}, which do not require a linear connection.

The paper is organized as follows. In Section 2, we recall some basic facts on Lie algebroids. Section 3 is devoted to the geometry of Lie algebroid prolongations and the constructions needed throughout the paper. In Section 4, we introduce semisprays, dynamical variations and the Jacobi
fields of \cite{Avendano-Kama-2026}. In Section~5 we characterize variations of
admissible curves with compatible transverse variations in terms of Lie
algebroid morphisms. In Section~6 we use this characterization to obtain
a variational description of the Jacobi sections of \cite{carinena2015jacobi} along
integral curves of a semispray.  In Section~7 we relate Jacobi sections and
Jacobi fields through the anchor map, and establish the conditions under
which a Jacobi field is the projection of a Jacobi section. Finally, in
Section~8 we present an example illustrating this correspondence.
%\cite{lee2018introduction,carinena2015jacobi} 

\section{Preliminaries on Lie algebroids}\label{sec2}

\begin{definition}
A Lie algebroid  is a triple $(p: A \rightarrow M, [ \cdot,\cdot ]_ A, \rho)$ where $p: A \rightarrow M$   is a vector bundle with a real Lie algebra bracket $[ \cdot,\cdot ]_A: \Gamma(A)\times \Gamma(A)\rightarrow \Gamma(A) $ defined on the $C^\infty(M)$-module of sections $\Gamma(A)$ and with a vector bundle homomorphism $\rho : A \rightarrow TM $, called anchor map, satisfying the Leibniz rule
\begin{equation*}
[\sigma, f \eta ]_A =  f[\sigma,\eta] + (L_{\rho(\sigma)}f)\eta,
\end{equation*}
for all $\sigma,\eta \in \Gamma(A)$, $f\in C^\infty(M).$
\end{definition}

Since the anchor map $\rho: A\rightarrow TM$ is a bundle morphism, it induces a morphism of $C^\infty(M)$-modules  $\rho: \Gamma(A) \rightarrow \mathfrak{X}(M)$. Moreover, an immediate consequence of this definition is that the anchor map  $\rho$ is a Lie algebra  homomorphism:
\begin{equation*}
\rho([ \sigma, \eta ]_A) = [\rho(\sigma),\rho(\eta) ]_{TM}, \ \ \ \  \sigma,\eta\in \Gamma(A).
\end{equation*}
The image of the anchor map  $\mathcal{D}:=\rho(A)$ is a smooth generalized distribution on $M$ which is integrable in the sense of Stefan-Sussmann \cite{stefan-1974,sussman-1973}. The integral submanifolds of  $\mathcal{D}$ define a singular foliation of $M$ called the \textit{characteristic foliation}.

If we take local coordinates $(U, x^1,\ldots , x^n)$ on $M$ and a local basis $\{ e_i \}_{i=1}^r$ of sections $\Gamma(U)$ we obtain local coordinates of $A$, called \textit{adapted coordinates}:
\begin{equation}
 A \ni a \mapsto (x^1,\ldots,x^n,y^1,\ldots,y^r),\label{adapcoord}
\end{equation} 
where $a = y^j \cdot e_j (x).$ The smooth functions $\rho^i_j$, $C^k_{ij}$ on $C^\infty(U)$ defined by
\begin{equation}
\rho(e_j)=\rho^i_j\frac{\partial}{\partial x^i}, \ \ \ \ [ e_i,e_j ]_A= C^k_{ij}e_k \label{strucfunc}
\end{equation} 
are called the structure functions of the Lie algebroid. These functions satisfy the so-called \textit{structure equations} on the Lie algebroid:
\begin{equation}
\rho^i_j\frac{\partial \rho^k_l }{\partial x^i}- \rho^i_l\frac{\partial \rho^k_j }{\partial x^i}= \rho_i^kC^i_{jl}, \ \ \ \sum_{\mathrm{cyclic}(j,l,s)}\left( \rho^i_j\frac{\partial C^k_{ls} }{\partial x^i} - C^k_{jt} C^t_{ls}\right)=0.\label{compconditions} 
\end{equation}

Now, we present some basic examples of Lie algebroids.

\begin{example}\label{tant-algebroid}
For any manifold $M$, the tangent bundle $\tau_M: TM \rightarrow M $ is a Lie algebroid with the standard Lie bracket of vector fields and anchor map given by the identity.
\end{example}

\begin{example}
Let $(M, \Pi)$ be a Poisson manifold, i.e., a manifold $M$ together with bivector $\Pi\in \mathfrak{X}^2(M)$  satisfying the Jacobi identity $$ [\Pi,\Pi ]_{\mathrm{SCH}}=0. $$ 
The Poisson bivector $\Pi $ induces the following vector bundle morphism

\begin{eqnarray*}
\Pi^\sharp: T^*M \rightarrow TM, \\
\Pi^\sharp(\omega):= \mathrm{i}_{\omega}\Pi.
\end{eqnarray*}
Moreover, it is well known that the following bracket
\begin{eqnarray*}
[\cdot,\cdot ]_\Pi&:& \Omega(M)\times \Omega(M) \rightarrow \Omega(M), \\
\left[ \sigma,\eta \right]_\Pi&:=& L_{\Pi^\sharp(\sigma)}\eta-L_{\Pi^\sharp(\eta)}\sigma-\mathrm{d}(\Pi(\sigma,\eta)),
\end{eqnarray*}
is a Lie bracket on the  $\mathbb{R} $-linear space of 1-forms $\Omega(M)$. Furthermore, the canonical cotangent bundle projection $\pi_M: T^*M \rightarrow M$ is a Lie algebroid with Lie bracket $[\cdot,\cdot ]_\Pi$ and anchor map $\Pi^\sharp: T^*M \rightarrow TM$. This example is discussed in full detail in \cite{CFM-2021}.
\end{example}

\begin{example}
An involutive distribution $\mathcal{D}\subset TM$ is a Lie algebroid with the inclusion $\mathcal{D} \hookrightarrow TM$ as anchor map and the usual bracket of vector fields as Lie bracket.
\end{example}

\begin{example}
Let  $\mathfrak{g}\ni \xi \mapsto \xi_M \in\mathfrak{X}(M) $ be a (right) action of the real Lie algebra $\mathfrak{g}$ on the manifold $M$. Its associated Lie algebroid consists of the trivial vector bundle $\mathrm{pr}_1:M \times\mathfrak{g} \rightarrow M$,
 with anchor map $\rho: M \times\mathfrak{g} \rightarrow TM $ defined by
$$\rho(m,\xi):=\xi_M(m). $$
 Identifying the sections of $\mathrm{pr}_1$ with the set of $\mathfrak{g}$-valued functions $v: M \rightarrow \mathfrak{g}$, one can define a Lie bracket by
 $$ [v,w](m):= [ v(m),w(m)] + \langle \mathrm{d}w, (v(m))_M \rangle (m) -    \langle \mathrm{d}v, (w(m))_M \rangle (m).$$
\end{example}

A section of the dual bundle $\pi: A^* \rightarrow M$ will be referred to as an $A$-form. More generally, elements of $\Gamma(\wedge^k A^*)$ will be called $A$-forms of degree $k$.

The set of 1-forms $\Gamma(A^*)$ can be identified with the space of linear functions $C_{\mathrm{lin}}^\infty(A) $ as follows: for every  $\theta \in \Gamma(A^*) $,
\begin{equation}
l_\theta(a):= \langle \theta_{p(a)},a \rangle. \label{linfunc} 
\end{equation}

The algebraic structures on a Lie algebroid allow us to define an exterior differential operator $\d_A: \Gamma(\wedge^{\bullet}A^*)\to\Gamma(\wedge^{\bullet}A^*)$, defined over the $C^\infty(M)$-module of $A$-forms $\Gamma(\wedge^\bullet A^*) $ as follows: if  $\omega\in\Gamma(\wedge^k A^*)$
\begin{eqnarray}
  \d_A\omega(\sigma_0,\ldots,\sigma_k) &:= & \sum_{i=0}^{k}(-1)^i\mathcal{L}_{\rho(\sigma_i)}(\omega(\sigma_0,\ldots,\hat{\sigma}_i,\ldots,\sigma_k))\nonumber \\
  && + \sum_{0\leqslant i<j \leqslant k}(-1)^{i+j}\omega([\sigma_i,\sigma_j]_A,\sigma_0,\ldots,\hat{\sigma}_i,\ldots,\hat{\sigma}_j,\ldots,\sigma_k),\label{diff-ext-operator}
\end{eqnarray}
where $\sigma_0,\ldots,\sigma_k\in\Gamma(A)$. The operator $\d_A$ is a graded derivation of degree one and a coboundary operator, i.e., $\d_A^2=0$,  of the cochain complex $(\Gamma(\wedge^{\bullet}A^*),\d_A)$.

%\begin{example}
%A Lie algebra $\mathfrak{g}$ can be thought of as a Lie algebroid $A \rightarrow M $ with $A= \mathfrak{g} $, the base manifold $M$ as a singleton $M= \{  \star \}, $ and anchor map the null map $\rho \equiv 0$.
%\end{example}

\begin{definition}
Given two Lie algebroids $(p_1: A_1 \rightarrow M_1, [ \cdot,\cdot ]_ {A_1}, \rho_1)$ and $(p_2: A_2 \rightarrow M_2, [ \cdot,\cdot ]_ {A_2}, \rho_2)$, a Lie algebroid morphism $\Phi: A_1 \rightarrow A_2$ over $\varphi: M_1 \rightarrow M_2$ is a vector bundle homomorphism $\Phi$ covering $\varphi$ such that the induce map on forms:
$$\Phi^* :\Gamma(\wedge^k A_2^*)\rightarrow \Gamma(\wedge^k A_1^*),$$
defined by $\langle \Phi^* (\sigma),a \rangle =\langle  \sigma_{\varphi(m)},\Phi(a) \rangle $, commutes with the exterior differential
$$ \mathrm{d}_{A_1}\circ \Phi^* = \Phi^* \circ\mathrm{d}_{A_2}.  $$ 
\end{definition}

The properties of Lie algebroid morphisms are summarized in the next proposition.

\begin{proposition}\label{cat-Lie-morph}
The set of Lie algebroids together with Lie algebroid morphisms form a category.
\end{proposition}

This is a standard fact which can be easily proven; see, for example, \cite{Kosmann-2008}.

\begin{example}
On any Lie algebroid  $(p: A \rightarrow M, [ \cdot,\cdot ]_ {A}, \rho)$, the anchor map $\rho: A \rightarrow TM$ is a Lie algebroid morphism between the Lie algebroid itself and the tangent bundle. From the formula \eqref{diff-ext-operator} for the differential exterior operator, one can deduce
$$ \rho^* \circ \mathrm{d} = \mathrm{d}_A \circ \rho^*. $$
\end{example}

%The next example of a morphism of Lie algebroids will be important, as a tool, in one of the main results in this paper.

\section{The geometry of Prolongation of a Lie algebroid}  

The prolongation of a Lie algebroid over its bundle projection turns out to be the natural analogue of the tangent bundle of a tangent bundle. It provides the natural geometric framework for the study of semisprays, second-order dynamics and variational objects on Lie algebroids. Here, we develop only the tools necessary for the paper. Further details on the Lie algebroid prolongation can be found in \cite{carinena2007reduction,deLeonEtAl2005,martinez2002lie,Martinez2001}.

For a given Lie algebroid $(p: A \rightarrow M, [ \cdot, \cdot ], \rho)$, its \emph{prolongation} over $A$  
is the vector bundle $p_A: \mathcal{T} A \rightarrow A   $ where
\begin{equation}
\mathcal{T} A :=\{ (a,b,v)\in A\times A \times TA \mid p(a)=p(b); \, \rho(b)=T_a p(v)   \}\label{prolalgebroid}
\end{equation}
and projection bundle given by % $p_A(a,b,v):=a$.
\begin{equation}
p_A(a,b,v):=a. \label{prolalgbundle}
\end{equation}
The fiber over $a\in A$ is isomorphic to the vector space
\begin{equation*}
(p_A)^{-1}(a) \simeq \{ (b,v) \in A_{p(a)}\times T_a A \mid  \rho(b)= T_a p(v) \}.
\end{equation*}

We also consider the map $\hat{p}:\mathcal{T} A \rightarrow A$ given by projection onto the second factor:
\begin{equation}
 \hat{p}(a,b,v):=b. \label{vert-proj} 
\end{equation}

%%%%. Vertical subbundle

%The \emph{vertical subbundle} of $\mathcal{T} A$ is 
%\begin{equation}\label{vert-subbundle}
%\mathrm{Vert}(\mathcal{T} A ) := \mathrm{ker}\ \hat{p} = \{ (a,b,v )\in \mathcal{T} A \mid b=0, \text{ and } \, T_ap(v)=0 \}.
%\end{equation}
%\begin{eqnarray}
%\mathrm{Vert}(\mathcal{T} A ) & := & \mathrm{ker}\ \hat{p}, \label{vert-subbundle} \\ 
%& =&\left\{ 
%\begin{array}{c|c}
%     \multirow{2}{*}{$(a,b,v )\in \mathcal{T} A$}& b=0,\nonumber  \\
%                                     & T_ap(v)=0. 
%\end{array}
%\right\}
%\end{eqnarray}
%Elements of $\mathrm{Vert}(\mathcal{T} A )$ are called \emph{vertical}. Moreover: %The mapping $p_A:\mathrm{Vert}(\mathcal{T} A ) \rightarrow A $ is a vector subbundle of $ \mathcal{T} A .$

%%%%%%%%%%%%%%%%%%%%%%%%%%

%%%%%%%%%%%%%%%%%%%%
%%%%%%%%%%%%%%%%%%%%
%%%%%% Vertical sections %%%%%
%%%%%%%%%%%%%%%%%%%%
%%%%%%%%%%%%%%%%%%%%

%\begin{definition} 
%A section $Z\in \Gamma( \mathcal{T} A )$ is called \emph{vertical} if  $Z(a)\in \mathrm{Vert}(\mathcal{T} A ) $ for all $a\in A$.
%\end{definition}
%Notice that a vertical section $Z$ in   $\mathcal{T} A  $ takes the form
%$$ Z = (\mathrm{id}_A, 0 , X),  $$
%with $X\in \mathfrak{X}(A)$ a vertical vector field on $A$, i.e., $Tp \circ  X=0.$

With the aid of the map $\hat{p}:\mathcal{T} A \rightarrow A$,  we can define a special kind of sections of $\mathcal{T}A$ called \emph{projectable} sections.

\begin{definition}
  A section $Z\in \Gamma( \mathcal{T} A )$ is said to be \emph{projectable} 
  if there exists $\sigma\in \Gamma(A)$ such that
\begin{equation}
\hat{p} \circ Z = \sigma \circ p.
\end{equation}
\end{definition}
One can deduce that every projectable section $Z$ has the form
$$ Z= (\mathrm{id}_A, \sigma\circ p , X), $$
where $\sigma \in \Gamma(A)$ and $X\in \mathfrak{X}(A)$ are such that $X$ and $\rho(\sigma)$ are $p$-related vector fields.
% is a vector field satisfying the condition
%$  Tp\circ X = \rho \circ \sigma \circ p.$ Here, $Tp:  TA \rightarrow TM $ denotes the tangent map of the projection bundle $p: A \rightarrow M$.

\paragraph{Lie algebroid structure of  $\mathcal{T}A$.} 
The key property that allows us to define a Lie algebroid structure on $\mathcal{T} A$ is the fact that 
every section of $\mathcal{T} A$ is locally written as a linear combination of projectable sections.
The Lie bracket can be easily defined in terms of projectable sections and extended by $\R$-linearity and the Leibniz rule.

To obtain a local basis of sections on $\mathcal{T} A$,  we have to consider the local coordinate system on the prolongation $\mathcal{T} A $ constructed as follows:
fix an adapted coordinate system $(x^i, y^k)$ on $A$ associated with a local chart $(U,x^1,x^2,\ldots,x^n)$ on $M$ and a basis $\{ e_l \}$ of local sections on $A$. 
Now, we choose a point $(a,b,v)\in \mathcal{T} A$ with $a,b \in p^{-1}(U)$. Since $p(a)=p(b)$, the points $a, b\in A$ have coordinates $(x^i,u^k) $ and $(x^i,y^j) $, respectively. Now, the condition $\rho(b)= T_a p(v) $ implies that $v$ has the form $\displaystyle v=y^j\rho^i_j(x)\left.\frac{\partial}{\partial x^i}\right|_a + v^l\left.\frac{\partial}{\partial y^l}\right|_a$. Thus, one can assign
\begin{equation}
\mathcal{T} A  \ni(a,b,v) \mapsto (x^i,u^k, y^j,  v^l),
\end{equation}
as coordinates. 
It follows from here that $\mathrm{rank}\,\mathcal{T} A  =2\,\mathrm{rank}\, A $ .
%and  $\mathrm{rank}\,(\mathrm{Vert}(\mathcal{T} A )) = \mathrm{rank}\, A$.

%%%%% Define vertical sections on the prolongations

Then, we define a local basis $\{\mathcal{E}_i, \Upsilon_j\}$ of sections of $\mathcal{T} A$ 
%associated with a local basis of section $\{ e_l \}$ on $A$ 
by formulas:
\begin{equation}\label{sectionprolon}
\mathcal{E}_j(a):= \left(a, e_j(p(a)),  \rho^i_j(x)\left.\frac{\partial}{\partial x^i}\right|_a  \right), \  \  \  \text{and}  \  \  \  \Upsilon_k(a):=\left(a, 0,\left.\frac{\partial}{\partial y^k}\right|_a\right)
\end{equation}
 These local sections are projectable:   $\hat{p} \circ \mathcal{E}_j = e_j \circ p$ and $\hat{p} \circ\Upsilon_k = 0$ .
%Moreover, the sections $\Upsilon_k$ are vertical and determine a local basis for the subspace of vertical sections.

Finally, the bracket of the set of projectable sections is presented in the next definition.

\begin{definition}
Let $Z_i=  (\mathrm{id}_A, \sigma_i \circ p , X_i) $ be projectable sections of $\mathcal{T} A $ with $\sigma_i\in\Gamma(A) $ for $i=1,2.$  
The Lie bracket of $Z_1$ and $Z_2$ is
\begin{equation}
[Z_1,Z_2](a):=(a,[\sigma_1,\sigma_2]\circ p(a),[X_1,X_2](a)).\label{brack-proj}
\end{equation}
\end{definition}

In particular, the Lie brackets of the elements of the basis \eqref{sectionprolon} are given by
\begin{equation}
[\mathcal{E}_i,\mathcal{E}_j ] = C^k_{ij}\mathcal{E}_k, \qquad  [\mathcal{E}_i, \Upsilon_k]=0, \quad \text{and} \quad  [\Upsilon_j, \Upsilon_k ]=0. \label{bracket-in-basis}
\end{equation}
%It follows from here, that the vertical bundle $\mathrm{Vert}(\mathcal{T} A ) $ is involutive
%$$ [\Gamma(\mathrm{Vert}(\mathcal{T} A ) ), \Gamma(\mathrm{Vert}(\mathcal{T} A ) )] \subset \Gamma(\mathrm{Vert}(\mathcal{T} A ) ). $$
The anchor map $\hat{\rho} : \mathcal{T}A \rightarrow TA$ is just the projection onto the third factor, 
\begin{equation}
\hat{\rho}(a,b,v) :=v,\label{anchor-prol}
\end{equation}
and can also be defined locally in terms of the basis \eqref{sectionprolon} by
\begin{equation}
\hat{\rho} (\mathcal{E}_j) = \rho^i_j\frac{\partial}{\partial x^i}, \qquad   \hat{\rho} (\Upsilon_k)=\frac{\partial}{\partial y^k}. \label{anchorinbasis}
\end{equation}

Therefore, the triple $(p_A: \mathcal{T} A \rightarrow A , [ \cdot, \cdot ]_{\mathcal{T} A}, \hat{\rho} \    )$, where the bracket $[ \cdot, \cdot ]_{\mathcal{T} A}$ is defined over projectable sections by \eqref{brack-proj} and the anchor map  $ \hat{\rho} $ given by \eqref{anchor-prol}, is a Lie algebroid called the prolongation Lie algebroid of $A$.

%%%%%%%%% morphism Lie algebroid property %%%%%
 
In the following result, we state an important property.% of the map 

\begin{proposition}\label{prop-phat}
The map $\hat{p}:\mathcal{T} A \rightarrow A$ given by \ref{vert-proj} is a Lie algebroid morphism over the vector bundle projection $p: A\rightarrow M$,  i.e., $\hat{p}$ is a vector bundle morphism satisfying
$$\hat{p}^{\,*} \circ \d_{\mathcal{T} A } = \d_A \circ \hat{p}^{\,*}.$$
\end{proposition}

The proof of this proposition can be found in \cite{Martinez2001}.

%Following \cite{Martinez2001},  the map $\hat{p}:\A:=\mathcal{T} A \rightarrow A$ is a Lie algebroid epimorphism over the projection $p:A\rightarrow M$.
%In consequence, we have that $\hat{p}^{\,*}$ is a cochain complex monomorphism: $\hat{p}^{\,*} \circ \d_{\A} = \d_A \circ \hat{p}^{\,*}$.
%Furthermore, since $\mathrm{Vert}(\A)=\ker\hat{p}$,  it follows that $\hat{p}^{\,*}:\wedge^{\bullet}A^*\to\wedge^{\bullet}\mathrm{Vert}(\A)^\circ$ is an isomorphism.

Lie algebroid morphisms can be lifted to the prolongations of Lie algebroids, as we briefly discuss. Assume that $(p_i: A_i \rightarrow M_i, [ \cdot,\cdot ]_ {A_i}, \rho_i)$, with $i=1,2$, are two Lie algebroids. Assume that $\Phi : A_1 \rightarrow A_2$ is a vector bundle morphism over $\varphi :M_1 \rightarrow M_2$.  Then,  the map $\mathcal{T}\Phi :\mathcal{T}A_1 \rightarrow \mathcal{T}A_2 $ given by
\begin{equation}
\mathcal{T}\Phi (a,b,v) = (\Phi(a),\Phi(b), T_a\Phi(v))\label{lift-phi}
\end{equation}
is also a vector bundle morphism over $\Phi$. Moreover, $\mathcal{T}\Phi$ inherits from $\Phi$ the property of being a Lie morphism, as we establish in the following proposition. The proof of this fact and more properties of the map $\mathcal{T}\Phi$ can be found in \cite{martinez-2008}.

\begin{proposition}\label{prop-lift-phi}
If $\Phi : A_1 \rightarrow A_2$ is a morphism of Lie algebroid then the map  $\mathcal{T}\Phi :\mathcal{T}A_1 \rightarrow \mathcal{T}A_2 $ \eqref{lift-phi} is a morphism of Lie algebroids,
\end{proposition}

In the prolongation of a Lie algebroid $ (p:A\rightarrow M, [ \cdot,\cdot], \rho) $, there exists map $\chi_A : \mathcal{T} A  \rightarrow \mathcal{T} A  $ with the property of $\chi_A  \circ \chi_A = \mathrm{Id} $, called the \textit{canonical involution} \cite{CortesEtAl2006,Martinez2001,deLeonEtAl2005,martinez-2008}. This map is given by
\begin{equation}
\chi_A (a,b,v) = (b,a, \bar{v}),\label{caninvol}
\end{equation}
where $\bar{v}\in T_bA$ such that $\rho(a)= T_bp(\bar{v})$ and
\begin{equation}
\bar{v}(l_\theta) = v(l_\theta) + \mathrm{d}_A \theta (a,b),\label{rel-invol}
\end{equation}
with $l_\theta$ the linear function \eqref{linfunc} defined for every section $\theta \in \Gamma(A^*).$

Given a section $\eta \in \Gamma(A)$, one can define the \textit{complete lift} of  $\eta $ as the section $\eta^c \in \Gamma(\mathcal{T} A)$ defined by

\begin{equation}
\eta^c(a) =  \chi_A(\eta\circ p(a), a, T_{p(a)} \eta (\rho(a))). \label{lift-section}
\end{equation}

\paragraph{The prolongation of the tangent Lie algebroid.} In this part, we discuss in some detail the prolongation of the tangent Lie algebroid of a smooth manifold. In particular, we are interested in showing how a vector field on the manifold induces a Lie algebroid morphism from the tangent Lie algebroid into its prolongation. This fact plays a crucial role in the proof of the main theorem in Section 6.

Let $M$ be a smooth manifold, and consider its tangent Lie algebroid $(\tau_M :  TM \rightarrow M, [ \cdot, \cdot]_{Lie}, \mathrm{id}_{TM}).$ The prolongation of the tangent Lie algebroid consists of the vector bundle $(\tau_M)^{TM} : \mathcal{T}  TM \rightarrow TM,$ where 
\begin{equation*}
\mathcal{T} TM=\Big\{(v,w,V)\in TM\times TM\times TTM \;\Big|\;
\tau_M(v)=\tau_M(w),\ \ T\tau_M(V)=w\Big\},
\end{equation*}
and $(\tau_M)^{TM}(v,w,V)=v. $ A projectable section on $\Gamma(\mathcal{T} TM)$ is of the form $Z(v)=(v, X(\tau_M(v), W(v))$, with $X\in \mathfrak{X}(M)$,
$W\in \mathfrak{X}(M)$ satisfying $T\tau_M W = X\circ \tau_M$. So, the bracket of two projectable sections $Z_1,Z_2\in\Gamma(\mathcal{T} TM)$ is given by
$$[Z_1,Z_2](v)=(v, [X_1,X_2]\circ \tau_M(v),[W_1,W_2]\circ (v)). $$
 Finally, the anchor map $\hat{\rho}_{\mathcal{T}TM}: \mathcal{T} TM \rightarrow TTM$ is the map $\hat{\rho}_{\mathcal{T}TM}(v,w,V)=V $. Observe that $\hat{\rho}_{\mathcal{T}TM}$ is an isomorphism of Lie algebroids with inverse mapping $(\hat{\rho}_{\mathcal{T}TM})^{-1}: TTM \rightarrow \mathcal{T} TM$  given by
 $$(\hat{\rho}_{\mathcal{T}TM})^{-1} (V)= \big(\tau_{TM}(V),\,T\tau_M(V),\,V\big). $$
 From the previous observations, we have the following consequences:
 \begin{itemize}
 \item For every section  $Z\in\Gamma((\mathcal{T} TM)^*)$ there exist a section $\varsigma \in \Omega(TM)$ such that $(\hat{\rho}_{\mathcal{T}TM})^*\varsigma = Z.$
 
 \item $\displaystyle (\hat{\rho}_{\mathcal{T}TM})^* \circ d_{TTM} =  d_{\mathcal{T} TM}\circ(\hat{\rho}_{\mathcal{T}TM})^*$.
 \end{itemize}
 
 Let $X\in \mathfrak{X}(M)$ be a vector field on $M$. For every $v \in TM $, the point $(X\circ \tau_M (v),v,T_v X (v))$ belongs to  $\mathcal{T} TM$. Thus, every vector field $X\in \mathfrak{X}(M)$ induces a map $\mathcal{C}_X:TM\longrightarrow\mathcal{T} TM,$ given by 
 
 \begin{equation}
 \mathcal{C}_X(v):=\Big(X(\tau_M(v)),\ v,\ T_{\tau_M(v)}X(v)\Big).\label{mapCX}
 \end{equation}
 By straightforward computations, one can show that  $ \mathcal{C}_X$ makes the following diagram 
 \[
\begin{tikzcd}[column sep=4em, row sep=3em]
TM \arrow[r, "\mathcal{C}_X"] \arrow[d, "\tau_M"'] & \mathcal{T} TM \arrow[d, "(\tau_{TM})^{TM}"] \\
M \arrow[r, "X"'] & TM
\end{tikzcd}
\]
commutative and $ \mathcal{C}_X$ is linear on each fiber. Furthermore, this mapping has the following property.
 
\begin{proposition}\label{prop-CX}
For every $X\in \mathfrak{X}(M)$, the map $\mathcal{C}_X$ given by \eqref{mapCX} is a Lie algebroid morphism over $X$, i.e., the pair $ (\mathcal{C}_X,X)$ is a vector bundle morphism such that
\begin{equation}
 d_{TM}\circ \mathcal{C}_X^{*}=\mathcal{C}_X^{*}\circ d_{\mathcal{T} TM}.\label{basirel}
\end{equation}
\end{proposition}
 \begin{proof}
 First, we show the compatibility of $ \mathcal{C}_X$ with the anchor maps. This means the commutativity of the diagram 
 \[
\begin{tikzcd}[column sep=4em, row sep=3em]
TM \arrow[r, "\mathcal{C}_X"] \arrow[d, "\rho_{TM}=\mathrm{id}"'] & \mathcal{T} TM \arrow[d, "\hat\rho_{\mathcal{T}TM}"] \\
TM \arrow[r, "TX"'] & TTM
\end{tikzcd}
\]
Indeed, for every $v\in TM$ we have
\begin{equation*}
\hat\rho_{\mathcal{T}TM}\big(\mathcal{C}_X(v)\big)
=\hat\rho_{\mathcal{T}TM}\Big(X(\tau_M(v)),\,v,\,T_{\tau_M(v)}X(v)\Big)
=T_{\tau_M(v)}X(v)
=TX\big(\rho_{TM}(v)\big).
\end{equation*}
Thus, we have $ \hat{\rho}_{\mathcal{T}TM}\circ \mathcal{C}_X = TX $. From here, we only need to show that relation \eqref{basirel} holds for 1-sections. We must recall that the vector field $X\in \mathfrak{X}(M)$ is a smooth map $X: M \rightarrow TM$ and then its tangent map $TX :  TM \rightarrow TTM$ is a morphism on the respective tangent Lie algebroids, which implies that $(TX)^* \circ d_{TTM}= d_{TM}\circ (TX)^*$. Now, let  $Z\in\Gamma((\mathcal{T} TM)^*)$ and $\varsigma \in \Omega(TM)$ such that $(\hat{\rho}_{\mathcal{T}TM})^*\varsigma = Z.$ By direct computation, we have
\begin{eqnarray*}
 d_{TM}\circ \mathcal{C}_X^{*}(Z)&=& d_{TM}\circ \mathcal{C}_X^{*}((\hat{\rho}_{\mathcal{T}TM})^*\varsigma) =d_{TM}\circ (\hat{\rho}_{\mathcal{T}TM}\circ \mathcal{C}_X)^{*}(\varsigma) =d_{TM}\circ (TX)^{*}(\varsigma),  \\
 &=& (TX)^* \circ d_{TTM}(\varsigma)= (\hat{\rho}_{\mathcal{T}TM}\circ \mathcal{C}_X)^{*} \circ d_{TTM}(\varsigma)= \mathcal{C}_X^{*}\circ(\hat{\rho}_{\mathcal{T}TM})^* \circ d_{TTM}(\varsigma),\\
 &=& \mathcal{C}_X^{*}\circ d_{\mathcal{T} TM}((\hat{\rho}_{\mathcal{T}TM})^*\varsigma)=\mathcal{C}_X^{*}\circ d_{\mathcal{T} TM}(Z).
\end{eqnarray*}
This completes the proof.
 \end{proof}

\section{Semisprays, dynamical variations and Jacobi fields}\label{sec2}

First, we introduce the fundamental notion of admissible curve on a Lie algebroid $ (p: A\rightarrow M, [ \cdot,\cdot]_A, \rho) $. 
  
 \begin{definition}
Let  $ (p:A\rightarrow M, [ \cdot,\cdot], \rho) $ be a Lie algebroid. A smooth curve $\vartheta: I \subset \mathbb{R} \rightarrow A$ is said to be \textit{admissible}  if
\begin{equation}
\frac{d \gamma}{dt}(t)=\rho (\vartheta(t)),
\end{equation}
where $\gamma = p\circ \vartheta$. The curve $\gamma$ is called the base curve of $\vartheta$.
\end{definition}
 In adapted coordinates  \eqref{adapcoord}, a curve $\vartheta(t)=(x^i(t),y^j(t))$  on $A$  is admissible if and only if $\dot{x}^i=\rho^i_j(x)y^j. $ In the literature these curves are also called $A$-paths, see \cite{CF-2003}.
 
It follows from the definition that the base curve of an admissible curve is completely contained in a leaf of the characteristic foliation. Moreover, every curve on $M$ contained in a leaf of the characteristic foliation is the base curve of some admissible curve, as we prove in the following lemma.

\begin{lemma}\label{existencecontlift}
Let $\mathcal{L}$ be a leaf of the characteristic foliation of a Lie algebroid  $ (p:A\rightarrow M, [ \cdot,\cdot], \rho) $ and let  $ \gamma :  I \rightarrow M$  be a smooth curve entirely contained in $\mathcal{L}$. Then, there exists an admissible curve $\vartheta : I \rightarrow A$ such that
\begin{equation*}
\rho(\vartheta(t))= \frac{d\gamma}{dt}.
\end{equation*} 
\end{lemma}
\begin{proof}
Let $A_\mathcal{L}=p^{-1}(\mathcal{L})$. Denote by $p_\mathcal{L}: A_\mathcal{L} \rightarrow \mathcal{L}$ the pull-back of $A$ via the immersion to $\iota_\mathcal{L}:  \mathcal{L} \hookrightarrow M $. Since the Lie bracket of section $\Gamma(A)$ induces a Lie bracket on the space of sections of $p_\mathcal{L}: A_\mathcal{L} \rightarrow \mathcal{L}$,  it has a transitive Lie algebroid structure with anchor map $\rho_\mathcal{L}$ induced by $\rho$, $\rho_\mathcal{L} := \rho\circ  \iota_\mathcal{L}$.  We also consider the Lie algebroid induced on $\mathcal{K}= \ker(\rho_{\mathcal{L}})$. Hence, we get an exact sequence of vector bundles over $\mathcal{L}$
\begin{equation}
0 \rightarrow \mathcal{K} \rightarrow A_\mathcal{L} \cong  T\mathcal{L}\oplus \mathrm{ker} \rho_\mathcal{L}  \rightarrow T\mathcal{L}\rightarrow 0. \label{shortexact}
\end{equation}
Because this exact sequence always splits, there exists a right splitting of $\rho_\mathcal{L}$, that is, a vector bundle map $\mathrm{h}: T\mathcal{L} \rightarrow A_\mathcal{L}$ such that $\rho_\mathcal{L}\circ \mathrm{h}= \mathrm{id}_{T\mathcal{L}} $. Now, we choose $\vartheta(t):= \mathrm{h}(\frac{d\gamma}{dt}(t))$. By construction, $p(\vartheta(t))=\gamma(t)$ and $\vartheta$ is an admissible curve,
 $$ \rho( \vartheta(t))= \rho_\mathcal{L}(\mathrm{h}(\frac{d\gamma}{dt}(t)))=\frac{d\gamma}{dt}(t). $$
\end{proof}

 \begin{example}
 Consider the tangent bundle $\tau_M :  TM \rightarrow M$ with its structure of Lie algebroid described in Example \ref{tant-algebroid}. Every smooth curve $\gamma: I \rightarrow M$ defines an admissible curve on $TM$ given by $\vartheta(t)= \dot{\gamma}(t).$
 \end{example}

Let $\vartheta: I \rightarrow A$ be an admissible curve on the Lie algebroid $ (p:A\rightarrow M, [ \cdot,\cdot], \rho) $. A \textit{variation}  of $\vartheta$ by admissible curves is a smooth map $\zeta: (-\varepsilon,\varepsilon)\times I \rightarrow A $, such that: 

\begin{itemize}

\item[(i)] $\zeta(0,t)=\vartheta(t)$,

\item[(ii)] for any $s\in (-\varepsilon,\varepsilon)$, the curve $t \mapsto \zeta(s,t) $ is admissible. 

\end{itemize}
 
The admissible curves are relevant objects on a Lie algebroid because any admissible curve is related to the velocity of its projection onto the base manifold by the anchor map. However, admissibility is a kinematic condition and, by itself, does not prescribe an evolution law. To describe dynamics on a Lie algebroid in terms of admissible curves, one may consider vector fields whose integral curves are admissible. This leads to the notion of semisprays on Lie algebroids, which play a role analogous to that of second-order differential equations on a tangent bundle.

\begin{definition}
A vector field $\mathcal{V}\in \mathfrak{X}(A)$ on the Lie algebroid  $ (p:A\rightarrow M, [ \cdot,\cdot], \rho) $ is called a \textit{semispray} if
\begin{equation}
Tp \circ \mathcal{V} = \rho, \label{defsemi}
\end{equation}
where $Tp: TA \rightarrow TM$ denotes the tangent map of the projection. A semispray  $\mathcal{V}\in \mathfrak{X}(A)$ is said to be a \textit{spray} if the vector field  $\mathcal{V}$ satisfies the homogeneous condition:
\begin{equation}
m^*_k  \mathcal{V}  = k\mathcal{V}, \ \ \ \ \text{for }k>0,\label{homogcond}
\end{equation}
where $m_k: A \rightarrow A$ is the dilation map given by $m_k(a):= ka$.
\end{definition}

\begin{remark}
Sprays on Lie algebroids, understood as vector fields on $A$ that lift the
anchor, also play a central role in \cite{cabrera2020local}, where a spray is used to
construct a local Lie groupoid integrating $A$. In that setting, the flow of
the spray gives rise to a path-exponential map with values in the groupoid,
while the exponential map considered in \cite{Avendano-Kama-2026} takes values in the
leaf of the base manifold. 
\end{remark}

We denote by $\phi^t_{\mathcal{V}}$ the local flow of a semispray $\mathcal{V}\in \mathfrak{X}(A)$. If $\vartheta: I \rightarrow A$ is an integral curve of  $\mathcal{V}$, its projection $\gamma(t)=p\circ \vartheta (t)$ onto $M$ will be called an \textit{integral curve in the base of } $\mathcal{V}$. A direct computation gives

\begin{equation*}
 \frac{d\gamma}{dt}(t)=\frac{d}{dt}(p\circ \vartheta(t)) = Tp \left(\frac{d \vartheta}{dt} \right)=  Tp (\mathcal{V}(\vartheta (t)) )=\rho(\vartheta (t)).
\end{equation*}
Thus, every integral curve of a semispray is admissible.

In adapted coordinates for $A$, a semispray $\mathcal{V}$  has the form
\begin{equation*}
\mathcal{V}(x,y)=y^j\rho^i_j(x)\frac{\partial}{\partial x^i}+\nu^k(x,y)\frac{\partial}{\partial y^k},
\end{equation*}
for some local functions $\nu^k(x,y) $ on $A$. If  $\mathcal{V}(x,y)$ is a spray, then $\nu^k(x,y)=\nu^k_{ln}(x)y^ly^n$. Note that every spray  $ \mathcal{V}$ vanishes at the zero section of the vector bundle $p: A \rightarrow M$, that is, $ \mathcal{V}(x,0)=0$. If $\vartheta(t)=(x^{i}(t), y^{k}(t))$ is an integral curve of $\mathcal{V}$ in local coordinates, then
\begin{eqnarray*}
\dot{x}^{i} = y^j\rho^i_j(x), \\
\dot{y}^{k} =\nu^k(x,y).
\end{eqnarray*}

\begin{remark}
The homogeneity condition  \eqref{homogcond} for a spray $\mathcal{V}$ can be expressed in terms of the  Liouville vector field $\mathbf{E}\in \mathfrak{X}(A)$, which locally has the form
\begin{equation*}
\mathbf{E}(x,y)=y^k\frac{\partial}{\partial y^k}.
\end{equation*}
Indeed, a semispray $\mathcal{V}\in \mathfrak{X}(A)$ is a spray if and only if 
\begin{equation*}
[ \mathbf{E},\mathcal{V} ]= \mathcal{V}.
\end{equation*}
\end{remark}

Now we fix an integral curve $\vartheta: I \rightarrow A$ of a semispray  $\mathcal{V}\in \mathfrak{X}(A)$ and the integral curve in the base $\gamma(t)=p\circ \vartheta (t)$. A \textit{dynamical variation}  of $\vartheta$ is a variation by integral curves of  $\mathcal{V}$, i.e. a smooth map $\zeta: (-\varepsilon,\varepsilon)\times I \rightarrow A $ such that: (i) $\zeta(0,t)=\vartheta(t)$, and (ii) for any $s\in (-\varepsilon,\varepsilon)$, the curve $\vartheta_s(t):=  \zeta(s,t) $ is an integral curve of  $\mathcal{V}$. In particular, $\zeta $ is a variation of $\vartheta(t)$ by admissible curves.

Then, we consider the smooth map $\alpha: (-\varepsilon,\varepsilon)\times I \rightarrow M $, given by

\begin{equation*}
\alpha := p\circ \zeta,
\end{equation*}
which is a variation by smooth curves of $\gamma(t)$. We will call $\alpha(s,t)$  a \textit{dynamical variation in the base}. Since $\zeta$ is a variation by integral curves of  $\mathcal{V}$, we have the following property:
\begin{equation}
\rho(\zeta(s,t))= \frac{\partial \alpha}{\partial t}(s,t).\label{poissrel}
\end{equation}

Let
\begin{equation}
Z(t)=\left.\frac{\partial \zeta}{\partial s}\right|_{s=0}, \ \ \text{and }\ \ \ J(t)=\left. \frac{\partial \alpha}{\partial s}\right|_{s=0}\label{variationalgamma}
\end{equation}
be the variational vector fields of variations $\zeta$  and  $\alpha$ respectively. In analogy with what happens in a Riemannian manifold $(M,g)$, one can say that the vector field 
\begin{equation*}
J(t)=\left. \frac{\partial \alpha}{\partial s}\right|_{s=0}
\end{equation*}
is a Jacobi field along the curve in the base $\gamma$, because it is the variational vector field of the dynamical variation in the base $\alpha(s,t).$ However,  such a characterization requires the existence of an explicit dynamical variation. It is therefore more convenient to formulate an equivalent definition directly in terms of vector fields along the corresponding integral curves. To state an equivalent notion but independent of dynamical variations, we need to look for a relation between the variational vector fields $Z(t)$ and $J(t)$. In \cite{Avendano-Kama-2026}, it is proved that the variational vector fields of the  dynamical variations $\zeta$ and $\alpha$, respectively,  are related by the equation
\begin{equation}
T\rho( Z(t)) = J^c(t), \label{eqbasic}
\end{equation}
where $T\rho: TA \rightarrow TTM $ denotes the tangent map of the anchor $\rho$ and $J^c(t)\in \Gamma( \left.TTM\right|_{\dot{\gamma}} )$ denotes the complete lift of the vector field $J(t)\in \Gamma( \left.TM\right|_{\gamma} )$.

Here, we give the definition of Jacobi fields for semisprays on Lie algebroids presented in \cite{Avendano-Kama-2026}.

\begin{definition}\label{defjacfield}
Let $\mathcal{V}\in \mathfrak{X}(A)$ be a semispray on the Lie algebroid $(p: A \rightarrow M, [ \cdot,\cdot ]_ A, \rho)$.  Fix an   integral curve in the base $\gamma: [0,T] \rightarrow M$ corresponding to an integral curve $\vartheta: [0, T]\subset \mathbb{R} \rightarrow A$ of  $\mathcal{V}$, i.e.,  $\gamma(t)=p(\vartheta(t))$. A vector field $J(t)$ along $\gamma$ is a \textit{Jacobi field} if there exists a vector field $Z(t)$ along $\vartheta(t)$ satisfying
\begin{eqnarray}
(L_\mathcal{V}Z)(t)&=&0, \label{JacCond1} \\
T\rho(Z(t)) &=& J^c(t). \label{JacCond2}
\end{eqnarray}
\end{definition}

\begin{remark}\label{remark2}
Condition \eqref{JacCond1} for a vector field $Z(t)$ along the integral curve $\vartheta(t)$ implies that $Z(t)$ is the variational vector field of some variation of $\vartheta(t)$ by integral curves of $\mathcal{V}$ which is equivalent to 
\begin{equation*}
T_a\phi_{\mathcal{V}}^t (\xi)=Z(t),
\end{equation*}
for some $\xi \in T_{\vartheta(0)}A$ and $T\phi_{\mathcal{V}}^t$ is the tangent map of $\phi_{\mathcal{V}}^t(\cdot) $ (see \cite{carinena2015jacobi}). Thus, we define Jacobi fields as variational vector fields of variations by integral curves in the base.
\end{remark}

\begin{remark}\label{remark3}
Condition \eqref{JacCond2} implies $ T p (Z(t)) = J(t)$, where $Tp: TA \rightarrow TM$ denotes the tangent map of $p: A \rightarrow M$.
\end{remark}

\begin{example}
Let $(M,g)$ be a Riemannian manifold and $\nabla$ the Levi-Civita connection of the metric $g$. On the tangent bundle $\tau_M : TM \rightarrow M$ with its canonical structure of Lie algebroid, 
the Riemannian metric induces the so-called geodesic spray $\mathcal{V}_g\in \mathfrak{X}(TM)$  \cite{chavel1995riemannian}, which is given locally by 
\begin{equation*}
\mathcal{V}_g(x,y) = y^i\frac{\partial}{\partial x^i} - \Gamma^k_{ij}(x)y^iy^j \frac{\partial}{\partial y^k}, 
\end{equation*}
where $\Gamma^k_{ij}$ are the Christoffel symbols of $\nabla$. The integral curves of $\mathcal{V}_g$ are the tangent lifts of geodesic curves on $M$. Hence, if $\zeta(s,t)$ is a dynamical variation of $\mathcal{V}_g$, then $\alpha(s,t)=p\circ\zeta(s,t)$ is a variation by geodesic curves. Therefore, the vector field
\begin{equation*}
J(t)= \frac{\partial \alpha}{\partial s}(0,t)
\end{equation*}
along $\gamma(t)=\alpha(0,t)$ is the classical variational characterization of Jacobi fields in Riemannian geometry.  As is well known, these variational vector fields are equivalently characterized as solutions of the classical Jacobi equation
 \begin{equation*}
 \nabla^2_{\dot{\gamma}} J + R(J,\dot{\gamma})\dot{\gamma}=0,
 \end{equation*}
 where $R$ denotes the curvature of the Levi-Civita connection, \cite{Lee-book-2018}.

On the other hand, one can verify that
$$  (L_{\mathcal{V}_g} J^c)(t)=0.$$
Since $\mathrm{id}_{TM}$ is the anchor map in the tangent Lie algebroid, the complete lift $J^c$ satisfies condition \eqref{JacCond2} automatically. Hence, the vector field $J(t)$ is also a Jacobi field along $\gamma(t)$ in the sense of Definition \ref{defjacfield}.
\end{example}

The previous example extends naturally from geodesic sprays to arbitrary semisprays on the tangent bundle.

\begin{example}
Consider the tangent Lie algebroid $(\tau_M : TM \rightarrow M, [ \cdot,\cdot ]_{\mathrm{Lie}}, \mathrm{id}_{TM}  ) $, where $[ \cdot,\cdot ]_{\mathrm{Lie}}$ stands for the canonical Lie bracket of vector fields on $M$. Let $\mathcal{V}\in \mathfrak{X}(TM)$ be a semispray on $TM$. In adapted coordinates, $\mathcal{V}$ takes the following form
\begin{equation*}
\mathcal{V}(x,y) = y^i\frac{\partial}{\partial x^i} + \nu^k(x,y) \frac{\partial}{\partial y^k}.
\end{equation*}
This kind of vector field is also known as a second-order differential equation on $TM$. Let $\vartheta : I \rightarrow TM$ be an integral curve of the semispray $\mathcal{V}$ and let $\gamma = p\circ \vartheta$. By Definition \ref{defjacfield}, a vector field $J(t)$ along $\gamma$ is a Jacobi field if and only if
\begin{equation*}
(L_{\mathcal{V}} J^c)(t)=0.
\end{equation*}
Therefore, in this case, we have recovered the definition of Jacobi fields for second-order differential equations on the tangent bundle given by Cari\~{n}ena, Gheorghiu and Mart\'inez in  \cite{carinena2015jacobi}, (Definition 2.6).
\end{example}

\section{Compatible transverse variations and Lie algebroid morphisms}

In the previous section, Jacobi fields for semisprays on Lie algebroids were introduced from a purely dynamical viewpoint as variational vector fields associated with dynamical variations of integral curves of a semispray. To connect this construction to the notion of Jacobi sections introduced by Cariñena, Gheorghiu and Martínez, we require additional geometric ingredients. This section introduces compatible transverse variations, which naturally complement dynamical variations and provide the appropriate framework for establishing the correspondence developed in the following sections.

Here, we consider variations  $\zeta: (-\varepsilon,\varepsilon)\times I \rightarrow A$ of admissible curves on $A$ with the property

\begin{itemize}

\item[(A)] the variation  $\alpha(s,t)=p(\zeta(s,t))$ of smooth curves on $M$ is contained in a fixed leaf $\mathcal{L}$ of the characteristic foliation.

\end{itemize}

\begin{definition}
Given a variation $\zeta$ by admissible curves satisfying the property (A), a transverse variation to $\zeta$ is a smooth map $ \varkappa : (-\varepsilon,\varepsilon)\times I \rightarrow A  $ such that: 
\begin{itemize}
\item[(i)] $p(\zeta(s,t))=p(\varkappa(s,t))$,
\item[(ii)] $\displaystyle \rho(\varkappa)=\frac{\partial \alpha}{\partial s}$. 
\end{itemize}
\end{definition}

The existence of a transverse variation to a variation $\zeta$ by admissible curves is ensured by any right splitting $\mathrm{h}: T\mathcal{L}\rightarrow A_{\mathcal{L}}$ of the short exact sequence of Lie algebroids  \eqref{shortexact}. In fact, the smooth mapping $ \varkappa(s,t) := \mathrm{h}(\frac{\partial \alpha}{\partial s})$ is a transverse variation of $\zeta(s,t)$. If $\rho$ is injective, there is a unique transverse variation to a given variation by admissible curves. On the contrary, a variation $\zeta$  admits many transverse variations to it. 

Let $\vartheta: I \rightarrow A$ be an admissible curve and $\gamma= p \circ \vartheta$. A smooth curve $\varsigma: I \rightarrow A$ such that $p\circ \varsigma = \gamma $ is called a section over $\gamma$. The set of sections over $\gamma$ will be denoted by $\Gamma_\gamma(A).$ Now, we consider the  map $\Xi_\vartheta :  \Gamma_\gamma(A) \rightarrow \Gamma( \left.TA\right|_{\vartheta})$ given by
\begin{equation}
\Xi_\vartheta(\varsigma)(t) = \hat{\rho}(\chi_A (\varsigma(t),\vartheta(t), \frac{d \varsigma}{dt}(t) ) ),\label{ximap}
\end{equation}
 where $\hat{\rho}: \mathcal{T}A \rightarrow TA $ denotes the anchor map of the prolongation \eqref{anchor-prol} and $\chi_A : \mathcal{T}A \rightarrow \mathcal{T}A $ is the canonical involution \eqref{caninvol}. 
 
 The complete lift of a section $\varsigma: I \rightarrow A$ over $\gamma$ is defined by
 \begin{equation}
 \varsigma^c(t) = \chi_A (\varsigma(t),\vartheta(t), \frac{d \varsigma}{dt}(t) ).
 \end{equation}
The map \eqref{ximap} can be written in terms of the complete lift by  $ \Xi_\vartheta(\varsigma)(t) = \hat{\rho}(  \varsigma^c(t) ).$  The map  \eqref{ximap} is extended to the case of a variation $\zeta: (-\varepsilon,\varepsilon)\times I \rightarrow A$ of admissible curves on $A$ and $\alpha = p\circ \zeta$.  A smooth map $c: (-\varepsilon,\varepsilon)\times I \rightarrow A$ such that $p\circ c = \alpha$ is called a section over $\alpha$ and the set of sections over $\alpha$ will be denoted by $\Gamma_\alpha(A).$ Then, we define the map $\Xi_\zeta :  \Gamma_\alpha(A) \rightarrow \Gamma( \left.TA\right|_{\zeta})$ by 
\begin{equation}
\Xi_\zeta(c) := \hat{\rho}(\chi_A (c,\zeta, \frac{\partial c}{\partial t} ) ).\label{ximapvar}
\end{equation}
Then, we introduce the \textit{tensor of variation} of a variation $\zeta: (-\varepsilon,\varepsilon)\times I \rightarrow A$ of admissible curves as the map $\Delta_\zeta : \Gamma_\alpha(A) \rightarrow \Gamma( \left.TA\right|_{\zeta})$ given by
\begin{equation}
\Delta_\zeta(\varkappa) :=  \frac{\partial \zeta}{\partial s} - \Xi_\zeta(\varkappa).\label{tensorvar}
\end{equation}

\begin{definition}
Given a variation $\zeta$ by admissible curves satisfying the property (A), a transverse variation to $\zeta$ is said to be compatible if
\begin{equation}
\Delta_\zeta(\varkappa) = 0.
\end{equation}
\end{definition}
 
The following result asserts the existence of a compatible transverse variation.

 \begin{proposition}\label{varzero} 
Given a variation $\zeta(s,t)$ by admissible curves with  property (A) and $\varkappa_0:(-\varepsilon,\varepsilon) \rightarrow A$ such that 
\begin{equation*}
\rho(\varkappa_0)(s)=\frac{\partial}{\partial s}(p\circ \zeta)(s,0),
\end{equation*}
there exists a unique transverse variation $\varkappa$ of $\zeta$ such that 
\begin{equation}
\Delta_\zeta(\varkappa)=0, \ \ \ \varkappa(s,0)=\varkappa_0(s).\label{uniquetrans}
\end{equation}
 
\end{proposition}
\begin{proof}
Suppose that the variation $\zeta$, the mapping $\alpha= p\circ \zeta$ and  the curve $\varkappa_0$ have the following local representations:
\begin{eqnarray}
\zeta(s,t)&=&\sum_{j=i}^ry^j(s,t)e_j(\alpha(s,t)), \  \  \  \  \varkappa_0(s)=\sum_{j=i}^rz_0^j(s)e_j\alpha(s,0),\nonumber \\
&&\alpha(s,t)=(x^1(s,t),\ldots, x^n(s,t)).\nonumber \\
\end{eqnarray}
By equation \eqref{tensorvar}, the variation $\varkappa(s,t)=\sum_{j=i}^rz^j(s,t)e_j\alpha(s,t)$ is compatible if and only if  for each $s\in(-\varepsilon,\varepsilon) $ the coefficient functions $z^j(s,t)$ satisfy the following time-dependent linear system of differential equations
\begin{eqnarray}
\frac{\partial z^k}{\partial t}   +y^iC^k_{ij}z^j & =&  \frac{\partial y^k}{\partial s}, \ \ \text{for all } k=1,\dots,r, \label{lineardiffeq}\\
z^k(s,0) &=& z_0^k(s).
\end{eqnarray}
Thus, the existence and uniqueness of $\varkappa$ follow from the existence and uniqueness theorem of solutions of linear differential systems. In order to prove that $\varkappa$ is a transverse variation of $\zeta$, we need to prove that
\begin{equation}
\frac{\partial x^k}{\partial s}(s,t)=z^j(s,t)\rho^k_j(\alpha(s,t)).\label{transcond} 
\end{equation}
By assumption,  $\frac{\partial x^k}{\partial t}=y^j(s,t)\rho^k_j(\alpha(s,t)).$ Then,
\begin{eqnarray*}
\frac{\partial}{\partial t}\left(\frac{\partial x^k}{\partial s}\right)&=&\frac{\partial}{\partial s}\left(\frac{\partial x^k}{\partial t}\right) = \frac{\partial y^j}{\partial s} \rho^k_j + y^j\frac{\partial \rho^k_j }{\partial x^l}\frac{\partial x^l}{\partial s}, \\
&&\text{ by equation  \eqref{lineardiffeq}}, \\
&=& \frac{\partial z^j}{\partial t}\rho^k_j    +y^lC^j_{lr}\rho^k_jz^r + y^j\frac{\partial \rho^k_j }{\partial x^l}\frac{\partial x^l}{\partial s}, \\
&&\text{ using the structure equations \eqref{compconditions}, } \\
&=& \frac{\partial}{\partial t}\left( z^j\rho^k_j \right) + y^j\frac{\partial \rho^k_j }{\partial x^l}\left( \frac{\partial x^l}{\partial s} - z^j\rho^k_l \right). 
\end{eqnarray*} 
Hence, for all $s\in (-\varepsilon,\varepsilon)$ we have the following time-dependend linear system of differential equations
\begin{equation*}
\frac{\partial}{\partial t}\left(\frac{\partial x^k}{\partial s}- z^j\rho^k_j \right) = y^j\frac{\partial \rho^k_j }{\partial x^l}\left( \frac{\partial x^l}{\partial s} - z^j\rho^k_l \right),\ \   \text{for all } k=1,\dots,n.
\end{equation*}
Since $\frac{\partial x^k}{\partial s}(s,0)- z^j(s,0)\rho^k_j(\alpha(s,0))=0 $, by uniqueness of solution of linear system of differential equations the relation \eqref{transcond} holds.
\end{proof}

Let $I, \tilde{I}  \subset \mathbb{R}$ be two open intervals. Consider the tangent bundle $\tau_{\tilde{I}\times I }: T\tilde{I}\times TI  \rightarrow \tilde{I}\times I  $ with its canonical Lie algebroid structure. Assume that $\Phi : T\tilde{I}\times TI   \rightarrow A$ is a vector bundle morphism over the smooth map $\alpha: \tilde{I}\times I   \rightarrow M$. Denote the coordinates $(s,t)$ in $\tilde{I}\times I  $ and let $\zeta(s,t)= \Phi(\left.\frac{\partial}{\partial t}\right|_{(s,t)})$ and  $\varkappa(s,t)= \Phi(\left.\frac{\partial}{\partial s}\right|_{(s,t)})$. Then, the vector bundle morphism can be written
\begin{equation}
\Phi(s,t) = \zeta(s,t) \mathrm{d}t + \varkappa(s,t)\mathrm{d}s. \label{morph-theta}
\end{equation}
Since $\Phi$ is vector bundle morphism over $\alpha$, we have that
$$ p\circ \zeta = p \circ \varkappa = \alpha. $$
The following question concerns the conditions under which the morphism $\Phi$ \eqref{morph-theta} is a morphism of Lie algebroids. In the next result, we present sufficient and necessary conditions under which $\Phi$ is a morphism of Lie algebroids. Thus, these conditions are related to variations of admissible curves with compatible transversal conditions.

\begin{theorem}\label{theo-morph-var}
The morphism $\Phi$ \eqref{morph-theta} is a morphism of Lie algebroids if and only if $\zeta(s,t)$ is a variation of admissible curves satisfying the property (A) and $\varkappa(s,t)$ is a compatible transverse variation of $\zeta$.
\end{theorem}
\begin{proof}
First, we assume that $\Phi$ is a Lie algebroid morphism. Then,  $\mathrm{d}( \alpha^* f )= \Phi^* \mathrm{d}_A f$,  for every $f\in C^\infty(M)$. A direct computation gives

$$ \langle \mathrm{d}f, \rho(\zeta)\rangle =\langle \mathrm{d}_Af, \zeta\rangle = \langle \Phi^* \mathrm{d}_A f, \frac{\partial}{\partial t} \rangle =   \langle \mathrm{d}(\alpha^* f ), \frac{\partial}{\partial t} \rangle = \langle \alpha^*\mathrm{d} f , \frac{\partial}{\partial t} \rangle =\langle \mathrm{d} f , \frac{\partial \alpha}{\partial t} \rangle .  $$
Thus, $\frac{\partial \alpha}{\partial t} = \rho(\zeta)$ and then, $\zeta(s,t)$ is a variation of admisible curves. By similar computations, we get that $\frac{\partial \alpha}{\partial s} = \rho(\varkappa)$, which implies that $\zeta(s,t)$ satisfies the property (A) and $\varkappa(s,t)$ is a transverse variation of $\zeta$. Now, for every $\theta\in \Gamma(A^*)$, we have $\mathrm{d} ( \Phi^* \theta) = \Phi^*( \mathrm{d}_A \theta)$ By one hand, we have
\begin{equation}
 \Phi^* (\mathrm{d}_A \theta )\left(\frac{\partial}{\partial s}, \frac{\partial}{\partial t}\right) = \mathrm{d}_A \theta (\varkappa, \zeta)\label{ec1}
\end{equation}
On the other hand,
\begin{eqnarray}
\mathrm{d}  (\Phi^* \theta)\left(\frac{\partial}{\partial s}, \frac{\partial}{\partial t}\right) &=& \frac{\partial}{\partial s} (l_\theta \circ \zeta ) - \frac{\partial}{\partial t} (l_\theta \circ \varkappa )=\langle \mathrm{d}(\zeta^*l_\theta), \frac{\partial}{\partial s} \rangle - \langle \mathrm{d}(\varkappa^*l_\theta), \frac{\partial}{\partial t} \rangle \nonumber\\
& =& \langle \zeta^*\mathrm{d}(l_\theta), \frac{\partial}{\partial s} \rangle - \langle \varkappa^* \mathrm{d}(l_\theta), \frac{\partial}{\partial t} \rangle = \frac{\partial \zeta}{\partial s}(l_\theta) - \frac{\partial \varkappa}{\partial t}(l_\theta) . \label{ec2}
\end{eqnarray}
Combining \eqref{ec1} and \ref{ec2}, we get
$$ \frac{\partial \zeta}{\partial s}(l_\theta) = \frac{\partial \varkappa}{\partial t}(l_\theta) +\mathrm{d}_A \theta (\varkappa, \zeta)$$
Taking into account \eqref {rel-invol}, we obtain
$$\chi_A\left( \varkappa, \zeta , \frac{\partial \varkappa}{ \partial t}\right) =  \left(\zeta ,  \varkappa, \frac{\partial \zeta}{\partial s}\right).  $$
By \eqref{anchor-prol} and \eqref{ximapvar}, $\displaystyle \Xi_\zeta(\varkappa)= \frac{\partial \zeta}{\partial s}$. Hence, $\Delta_\zeta(\varkappa)=0 $ and then $\varkappa$ is a compatible transverse variation os $\zeta$. The converse of this result also holds because all the previous arguments are reversible.
\end{proof}

Under the assumptions of Theorem \ref{theo-morph-var}, we assume that $\zeta$ is a variation of the admissible curve $\vartheta : I \rightarrow A$, i.e., $\zeta(0,t)=\vartheta(t)$.
Let  $\kappa \in \Gamma_\gamma(A)$, with $\gamma = p\circ \vartheta$, given by $\kappa(t)= \varkappa(0,t)$. Then, the infinitesimal variation $Z(t)$ of the variations $\zeta(s,t)$ is given by
\begin{equation}
Z(t)=\frac{\partial \zeta}{\partial s}(0,t)=\Xi_\vartheta(\kappa)(t)= \hat{\rho}(\kappa^c_\vartheta(t)).
\end{equation} 

The previous theorem shows that compatible transverse variations admit an intrinsic characterization in terms of Lie algebroid morphisms. This characterization will allow us to reinterpret the variational definition of Jacobi sections introduced by Cariñena, Gheorghiu and Martínez entirely in terms of compatible transverse variations.

%\begin{remark}
%Theorem \ref{theo-morph-var} is closely related to Proposition 2 in Martínez \cite{martinez-2008}, where Lie algebroid morphisms are characterized through the canonical involution $\chi_A$. The proof of our theorem is entirely intrinsic and based on the variation tensor, which provides a direct geometric interpretation of the compatibility condition.
%\end{remark}

\begin{remark}
The idea that a pair $(\zeta,\varkappa)$ defining a homotopy of $A$-paths should be
viewed as a morphism $TI\times TI\to A$ already appears in
\cite[\S 1.3]{CF-2003}, and a characterization of Lie algebroid morphisms in
terms of the canonical involution is given in \cite[Prop.~2]{martinez-2008}.
Theorem~\ref{theo-morph-var} makes this precise in the setting needed here:
it gives a coordinate-free proof that a vector bundle morphism
$\Phi=\zeta\,dt+\kappa\,ds$ is a Lie algebroid morphism if and only if
$\zeta$ is a variation by admissible curves and $\kappa$ is a compatible
transverse variation of $\zeta$, with compatibility expressed through the
variation tensor $\Delta_\zeta$.
\end{remark}

One can generate variations of admissible curves with compatible transverse variation by choosing a section $\sigma \in \Gamma(A)$ and a smooth curve $\nu: \widetilde{I} \rightarrow \mathcal{L} \subset M$, where $L$ is a leaf of the characteristic foliation of $\rho$. Following Cari\~{n}ena, Gheorghiu and Mart\'inez in \cite{carinena2015jacobi}, we say that the flow of a section $\sigma$ is the pair $(\phi_t,\varphi_t) $ where $\varphi_t$ denotes the local flow fo the vector field $\rho(\sigma)$ and $\phi_t$ the flow of $\hat{\rho}(\sigma^c)$. Here, $\sigma^c$ is the complete lift of $\sigma$ defined in \eqref{lift-section}. The flow  $(\phi_t,\varphi_t) $ of a section $\sigma\in\Gamma(A)$ has the following properties:

\begin{itemize}
\item The pair  $(\phi_t,\varphi_t) $ is a morphism of Lie algebroids.

\item   $ \sigma\circ \phi_t= \varphi_t \circ \sigma$.

\item The vector fields $\rho(\sigma)$ and $\hat{\rho}(\sigma^c)$ are $p$-related. 

\end{itemize} 

We define $\zeta : \widetilde{I}\times I \rightarrow A$ by
\begin{equation}
\zeta(s,t):= \sigma(\varphi_t(\nu(s))),\label{var-sec}
\end{equation}
where $I$ is the common interval of definition of the flow $\varphi_t$ throughout the points $\nu(s)$ with $s\in \widetilde{I}.$ Observe that $\alpha(s,t) = p\circ\zeta(s,t)=\varphi_t(\nu(s))$. Then
\begin{equation*}
\frac{\partial \alpha}{\partial t}(s,t)= \rho(\sigma(\varphi_t(\nu(s)))  ) = \rho(\zeta(s,t)).
\end{equation*}
Therefore, $\zeta(s,t)= \sigma(\varphi_t(\nu(s))) $ is a variation of admissible curves satisfying the property (A). By Lemma \ref{existencecontlift}, there exists a smooth curve $\varkappa_0: \widetilde{I} \rightarrow A$ such that $p\circ \varkappa_0 = \nu $ and $ \displaystyle \frac{d \nu}{ds} = \rho(\varkappa_0)$. By Proposition \ref{varzero} there exists a compatible transverse variation $\varkappa : \widetilde{I}\times I \rightarrow A$ such that $\varkappa(s,0)=\varkappa_0(s).$

The following result states an explicit formula for a compatible transverse variation of $\zeta$.

\begin{lemma}
The map $\varkappa:  \widetilde{I}\times I \rightarrow A $ given by 
\begin{equation*}
\varkappa (s,t)= \phi_t(\varkappa_0(s)), 
\end{equation*}
is a compatible transverse variation of $\zeta$ \eqref{var-sec}.
\end{lemma}
\begin{proof}
First, we shall prove that $\varkappa$ is a transverse variation of $\zeta$. Using that $\phi_t$ is a morphism of vector bundles over $\varphi_t$, we get
\begin{equation}
p(\varkappa (s,t)) = p(\phi_t(\varkappa_0(s)))= \varphi_t(p(\varkappa_0(s)))=\varphi_t(\nu(s))=p(\zeta(s,t)).
\end{equation}
Since $ (\phi_t,\varphi_t)  $ is a morphim of Lie algebroids,  $ \phi^*_t \mathrm{d}_A f= \mathrm{d}_A \circ \varphi_t^* f  $,  for all $f\in C^\infty(M)$. Then, 
\begin{eqnarray*}
\langle \mathrm{d}f,\rho(\varkappa (s,t)) \rangle&=& \langle \mathrm{d}_Af \phi_t(\varkappa_0(s)) \rangle =  \langle \phi^*_t \mathrm{d}_A f, \varkappa_0(s)\rangle = \langle \mathrm{d}_A \circ \varphi_t^* f  , \varkappa_0(s)\rangle,\\
&=& \langle \mathrm{d} \circ \varphi_t^* f  ,\rho( \varkappa_0(s))\rangle = \langle \varphi_t^*  \circ  \mathrm{d}  f,\frac{d\nu}{ds}(s)\rangle = \langle \mathrm{d} f  ,\frac{\partial}{\partial s}(\varphi_t(\nu( s))\rangle.
\end{eqnarray*}
Due to $f$ is arbitrary, $\displaystyle \frac{\partial}{\partial s}(\varphi_t(\nu( s))= \rho(\varkappa (s,t))$ which implies that $\varkappa$ is a transverse variation of $\zeta$. Now, we shall prove that $\varkappa$ is compatible. By direct computation, we get
\begin{eqnarray*}
\frac{\partial \varkappa}{\partial t} &= &\hat{\rho}(\sigma^c(\varkappa(s,t))),\\
&&\ \text{by \eqref{lift-section} we have,} \\
 &=& \hat{\rho}\left( \chi_A\left( \zeta(s,t),\varkappa(s,t),\frac{\partial \zeta}{\partial s}(s,t) \right)  \right).
\end{eqnarray*}
This equality implies
$$ \chi_A\left( \zeta(s,t),\varkappa(s,t),\frac{\partial \zeta}{\partial s}(s,t) \right) =  \left( \varkappa(s,t),\zeta(s,t),\frac{\partial \varkappa}{\partial t}(s,t) \right).$$
Since $\chi_A\circ \chi_A = \mathrm{id}, $ we obtain
$$ \chi_A\left( \varkappa(s,t),\zeta(s,t),\frac{\partial \varkappa}{\partial t}(s,t) \right)=\left( \zeta(s,t),\varkappa(s,t),\frac{\partial \zeta}{\partial s}(s,t) \right)  .$$
It follows from here that $\displaystyle \frac{\partial \zeta}{\partial s}(s,t) = \Xi_\zeta (\varkappa)$ which is equivalte to $\Delta_\zeta(\varkappa)=0$. This proves that $\varkappa$ is compatible.
\end{proof}

Furthermore, if $\Phi : T\tilde{I}\times TI   \rightarrow A$ is a morphism of Lie algebroids over the map $\alpha: \tilde{I}\times I   \rightarrow M$ if given by 
\begin{equation}
\Phi(s,t) = \sigma(\alpha (s,t)) \mathrm{d}t + \varkappa(s,t)\mathrm{d}s. \label{morph-section}
\end{equation}
with $\sigma \in \Gamma(A)$, then the map $\alpha(s,t)$ represents a variation of integral curves of the vector field $\rho(\sigma)$,
$$ \frac{\partial \alpha}{\partial t} (s,t)= \rho( \sigma)(\alpha(s,t)).  $$
That is, we can write $\alpha(s,t)= \varphi_t(\nu(s) )$ for some smooth curve $\nu: \tilde{I} \rightarrow L \subset M$. 

The previous discussions motivate the following definition, which is also presented in \cite{carinena2015jacobi}.

\begin{definition}
A 1-parameter family of integral curves of a section $\sigma\in \Gamma(A)$ is a Lie algebroid morphism $\Phi : T\tilde{I}\times TI   \rightarrow A$  over  $\alpha: \tilde{I}\times I   \rightarrow M$ of the form
\begin{equation}
\Phi(s,t) = \sigma(\alpha (s,t)) \mathrm{d}t + \varkappa(s,t)\mathrm{d}s. \label{1-para-int}
\end{equation}
The section $\kappa: I \rightarrow A$ given by $\kappa(t)=\varkappa(0,t)$ is called the infinitesimal variation of the 1-parameter family $\Phi$.
\end{definition}

\section{Jacobi sections for SODEs on Lie algebroids}

Once we have established the geometric framework of compatible transverse variations, we are ready to revisit the notion of a Jacobi section introduced by Cariñena, Gheorghiu and Mart\'inez \cite{carinena2015jacobi}.

Let $(p: A \rightarrow M, [ \cdot,\cdot ]_ A, \rho)$ and consider the prolongation of the Lie algebroid $p_A: \mathcal{T} A \rightarrow A $ given by  \eqref{prolalgebroid} and \eqref{prolalgbundle}.

\begin{definition} \label{def-sode}
A  second-order differential equation (SODE) of the Lie algebroid $A$ is a section $S\in \Gamma(\mathcal{T}A)$ in the prolongation given by 
\begin{equation}
S(a)=(a,a,\mathcal{V}(a)).\label{SODE}
\end{equation}
By an integral curve of a second-order differential equation $S$ is an integral curve of the vector field $\hat{\rho}(S)=\mathcal{V} \in \mathfrak{X}(A).$
\end{definition}

There is a close relation between semisprays and second-order differential equations on a Lie algebroid. It follows from the definition of the prolongation of a Lie algebroid that if $S\in \Gamma(\mathcal{T}A)$ is a SODE, then, the vector field  $\hat{\rho}(S)=\mathcal{V} $ is a semispray on $A$, i.e., $Tp \circ \mathcal{V} = \rho.$ Conversely, every semispray $\mathcal{V} \in \mathfrak{X}(A)$ induces a SODE by formula \eqref{SODE}. Hence, an integral curve of a SODE on $A$ is precisely an integral curve of a semispray on $A$.

\begin{remark}
An equivalent definition to Definition \ref{def-sode} is given in \cite{Martinez2001} in terms of the vertical endomorphism on the prolongation and the Euler section. This definition is analogous to the one given for the notion of second-order differential equation on the tangent bundle of a smooth manifold. Throughout this paper we adopt Definition \ref{def-sode} because it establishes the correspondence between semisprays on Lie algebroids and second-order differential equations on the prolongation that will be used below.
\end{remark}

Let $S\in \Gamma(\mathcal{T}A)$ be a SODE on $A$ and suppose that $ \Theta: \widetilde{I}\times I \rightarrow \mathcal{T}A$ is a 1-parameter family of integral curves of $S$. Now, we recall the notion of Jacobi section for second-order differential equations due to Cariñena, Gheorghiu and Mart\'inez, which was formulated in terms of infinitesimal variations of 1-parameter families of integral curves.

\begin{definition}\label{jacobi-sec-def}
Given an integral curve $\vartheta: I \rightarrow A$ of a SODE $S$ on the Lie algebroid $A$, a section $\kappa$ of $A$ along $\gamma= p \circ \vartheta$ is a Jacobi section along $\vartheta$ if its complete lift $\kappa^c_\vartheta(t)$ of $\kappa$ is the infinitesimal variation of $\Theta: \widetilde{I}\times I \rightarrow \mathcal{T}A$, a 1-parameter family of integral curves of $S$.
\end{definition}

%\begin{remark}
%Cariñena, Gheorghiu and Martínez also characterized Jacobi sections as solutions of a Jacobi-type equation on the prolongation of the Lie algebroid. In the present work, however, we focus on their variational characterization, since this point of view extends naturally to arbitrary semisprays via compatible transverse variations.
%\end{remark}

\begin{remark}
Cari\~{n}ena, Gheorghiu and Mart\'inez  also characterized Jacobi sections as solutions of a second-order equation, $\nabla\nabla W+\Phi(W)=0$, where the dynamical covariant derivative $\nabla$ and the Jacobi endomorphism $\Phi$ of the semispray are defined by means of the prolongation \cite[Thm.~3.7]{carinena2015jacobi}. In this paper we focus on their variational characterization, which we reformulate in terms of variations of integral curves on $A$.
\end{remark}

One of the purposes of this paper is to define Jacobi sections in an equivalent form to Definition \ref{jacobi-sec-def} but in terms of dynamical variations of integral curves and compatible transverse variations.  To achieve this goal, we provide a geometric interpretation of 1-parameter families of integral curves of SODEs.

If  $\Theta: T (\widetilde{I}\times I) \rightarrow \mathcal{T}A$ is a 1-parameter family of integral curves of $S$, then, $\Theta$ has the following form

\begin{equation}
\Theta = S(\zeta(s,t))dt + B(s,t)ds. \label{1-parameter-fam}
\end{equation}
Then,  $\Theta: T(\widetilde{I}\times I) \rightarrow \mathcal{T}A$ is a morphism of Lie algebroids over the map $\zeta: \widetilde{I}\times I \rightarrow A$. Besides, $\zeta$ is a variation by integral curves of the semisprays $\mathcal{V}=\hat{\rho}(S)$ or a dynamical variation of the integral curve $\vartheta(t):=\zeta(0,t)$. The coefficient $B(s,t)$ in \eqref{1-parameter-fam} has the following form
\begin{equation}
B(s,t)=\left(\zeta(s,t),\varkappa(s,t),\frac{\partial \zeta}{\partial s}(s,t) \right).\label{B-coef}
\end{equation}
This easily follows from the following identities: $\displaystyle p_A(B(s,t))= p_A(S(\zeta(s,t)))=\zeta(s,t) $ and   $\displaystyle\hat{\rho}(B(s,t))=\Theta\left(\frac{\partial}{\partial s}\right)= \frac{\partial \zeta}{\partial s}$. Moreover, since, $B(s,t)\in \mathcal{T}A$, we also have $ p(\varkappa(s,t))=p(\zeta(s,t))=\alpha(s,t)$ and $\rho(\varkappa)= Tp(\frac{\partial \zeta}{\partial s}(s,t)) =\frac{\partial \alpha}{\partial s}(s,t)$. This means that $\varkappa$ is a transverse variation of the dynamical variation $\zeta$ of $\vartheta$. Then, according to the Definition \ref{jacobi-sec-def}, the section $\kappa(t)=\varkappa(0,t)$ along $\gamma=p\circ \vartheta$ is a Jacobi section because of
\begin{equation*}
B(0,t)=(\vartheta(t),\kappa(t),Z(t))=\chi_A\left(\varkappa,\vartheta,\frac{d \kappa}{dt}(t) \right) =\kappa^c_\vartheta(t).
\end{equation*}
The previous discussions give rise to the relations between 1-parameter families of integral curves of a SODE and Jacobi sections with dynamical variations and their compatible transverse variations. 

In fact, we have the following result.

\begin{theorem}\label{theo-jacobi-transv}
Let $S\in \Gamma(\mathcal{T}A)$ be a SODE on $A$ and $\vartheta: I \rightarrow A$ an integral curve of  $S$. A section $\kappa$ of $A$ along $\gamma= p \circ \vartheta$ is a Jacobi section along $\vartheta$ if and only if there exists a dynamical variation $\zeta: \widetilde{I}\times I \rightarrow A$ of $\vartheta$ satisfying property (A) and a compatible transverse variation $\varkappa: \widetilde{I}\times I \rightarrow A$ of $\zeta$ such that $\varkappa(0,t)=\kappa(t)$.
\end{theorem}
\begin{proof}
First, we assume that $\kappa$ is a Jacobi section along $\vartheta$. Then, the complete lift of $\kappa(t)$ is the infinitesimal variation of $ \Theta: T( \widetilde{I}\times I) \rightarrow \mathcal{T}A$, a 1-parameter family of integral curves of $ S$. Assume that $\Theta$ is given by \eqref{1-parameter-fam} with coefficient $B$ of the form \eqref{B-coef}. %We have proven that $\zeta$ is a dynamical variation of $\vartheta$ with $\varkappa$ a transversal variarion of $\zeta$ such that $\kappa(t)=\varkappa(0,t)$. We only need to prove that $\varkappa$ is compatible, i.e., $\Delta_\zeta(\varkappa) = 0$. 
Recall the map $\hat{p}:\mathcal{T} A \rightarrow A$ given by \eqref{vert-proj} and define $\Phi: \widetilde{I}\times I \rightarrow A$ by $\Phi := \hat{p}\circ \Theta$. By Proposition \ref{cat-Lie-morph} and Proposition \ref{prop-phat}, $\Phi$ is a Lie algebroid morphism over $\alpha = p\circ \zeta$. A direct computation gives, $ \Phi\left( \frac{\partial}{\partial t} \right) = \hat{p}(S(\zeta(s,t)))= \zeta(s,t)$ and $ \Phi\left( \frac{\partial}{\partial s} \right) = \hat{p}(B(s,t))= \varkappa(s,t)$. Hence, 
\begin{equation}
 \Phi= \zeta(s,t) dt + \varkappa(s,t) ds.\label{theta-aux}
 \end{equation}
By Theorem \ref{theo-morph-var}, $\zeta: \widetilde{I}\times I \rightarrow A$ is a variation by admissible curves satisfying property (A) and $\varkappa: \widetilde{I}\times I \rightarrow A$ is a compatible transverse variation such that $\varkappa(0,t)=\kappa(t)$. Moreover, $\Theta(\partial/\partial t)=S(\zeta(s,t))$ and $\hat\rho$ maps it
to $\partial\zeta/\partial t$, so $\partial\zeta/\partial t=V(\zeta)$. Hence, for each $s\in \widetilde{I}$, the curve
$\zeta(s,\cdot)$ is an integral curve of $V=\hat\rho(S)$, which implies that $\zeta: \widetilde{I}\times I \rightarrow A$ is a dynamical variation of $\vartheta$.

%%%%%%%%%%
Now, we shall prove the converse. Assume that $\zeta: \widetilde{I}\times I \rightarrow A$ is a dynamical variation of $\vartheta$ satisfying property (A) and $\varkappa: \widetilde{I}\times I \rightarrow A$ is a compatible transverse variation such that $\varkappa(0,t)=\kappa(t)$, and construct the maps $\Phi$ by \eqref{theta-aux}. It follows from Theorem  \ref{theo-morph-var} that $ \Phi $ is a Lie algebroid morphism. Consider the mapping $\mathcal{T}\Phi:  \mathcal{T} T(\tilde{I}\times I) \rightarrow \mathcal{T} A$ defined over the corresponding prolongations of Lie algebroids by formula  \eqref{lift-phi}. By Proposition \ref{prop-lift-phi}, $\mathcal{T}\Phi$ is also a Lie algebroid morphism. Then, we choose the vector field $X=\frac{\partial}{\partial t} \in \mathfrak{X}(\tilde{I}\times I)$ and consider the map $\mathcal{C}_{\frac{\partial}{\partial t}}: T(\tilde{I}\times I) \rightarrow \mathcal{T} T(\tilde{I}\times I)$ by means of formula \eqref{mapCX}. $\mathcal{C}_{\frac{\partial}{\partial t}}$ is a Lie algebroid morphism by Proposition \ref{prop-CX}. Thus, the map $ \Theta :  T(\tilde{I}\times I) \rightarrow \mathcal{T} A$ given by
$$ \Theta :=  \mathcal{T}\Phi\circ \mathcal{C}_{\frac{\partial}{\partial t}}  $$
is also a Lie algebroid morphism. Evaluating the morphism $\Theta$ at the section $\frac{\partial}{\partial t}$, we get 
\begin{align*}
\Theta\!\left(\frac{\partial}{\partial t}\right)&=\mathcal T\Phi\circ C_{\partial/\partial t}\!\left(\frac{\partial}{\partial t}\right)
%&=\left(\Phi\!\left(\frac{\partial}{\partial t}\right),\ \Phi\!\left(\frac{\partial}{\partial t}\right),\
%T\Phi\!\left(T\!\left(\frac{\partial}{\partial t}\right)\!\left(\frac{\partial}{\partial t}\right)\right)\right)\\
=\left(\zeta(s,t),\ \zeta(s,t),\
T\big(\Phi(\partial/\partial t)\big)\!\left(\frac{\partial}{\partial t}\right)\right)\\
&=\left(\zeta(s,t),\ \zeta(s,t),\ \frac{\partial\zeta}{\partial t}\right)
=S(\zeta(s,t)).
\end{align*}
Analogously, evaluating the morphism $\Theta$ at the section $\frac{\partial}{\partial s}$, we obtain 
\begin{equation*}
\Theta\!\left(\frac{\partial}{\partial s}\right)
=\left(\Phi\!\left(\frac{\partial}{\partial t}\right),\ \Phi\!\left(\frac{\partial}{\partial s}\right),\
T\zeta\!\left(\frac{\partial}{\partial s}\right)\right)
=\left(\zeta(s,t),\ \varkappa(s,t),\ \frac{\partial\zeta}{\partial s}\right)
=B(s,t).
\end{equation*}
This implies that $\Theta$ is a 1-parameter family of integral curves of $S$. Then, $\kappa(t)$ is a Jacobi section along $\vartheta$ because 

\begin{equation*}
B(0,t)=(\vartheta(t),\kappa(t),Z(t))=\chi_A\left(\varkappa,\vartheta,\frac{d \kappa}{dt}(t) \right) =\kappa^c_\vartheta(t).
\end{equation*}

%Following \cite{Martinez2001},  the map $\hat{p}:\A:=\mathcal{T} A \rightarrow A$ is a Lie algebroid epimorphism over the projection $p:A\rightarrow M$.
%In consequence, we have that $\hat{p}^{\,*}$ is a cochain complex monomorphism: $\hat{p}^{\,*} \circ \d_{\A} = \d_A \circ \hat{p}^{\,*}$.
%Furthermore, since $\mathrm{Vert}(\A)=\ker\hat{p}$,  it follows that $\hat{p}^{\,*}:\wedge^{\bullet}A^*\to\wedge^{\bullet}\mathrm{Vert}(\A)^\circ$ is an isomorphism.

\end{proof}

From Theorem \ref{theo-jacobi-transv}, one can conclude that dynamical variations with compatible transverse variations provide the natural variational mechanism underlying Jacobi sections associated with semisprays on Lie algebroids.  With this result, dynamical variations and compatible transverse variations become the fundamental geometric objects, while the prolongation appears naturally as the geometric framework supporting the Definition \ref{jacobi-sec-def}, given originally in \cite{carinena2015jacobi} by Cariñena et al.

In consequence, the original definition of Cariñena, Gheorghiu and Martínez admits an equivalent formulation, expressed entirely in terms of dynamical variations and compatible transverse variations. Thus, in view of Theorem  \ref{theo-jacobi-transv}, Jacobi sections for semisprays (or SODEs) on Lie algebroids can equivalently be defined as follows.

\begin{definition}\label{defnew-jacobi}
Given an integral curve $\vartheta: I \rightarrow A$ of a SODE $S$ on the Lie algebroid $A$, a section $\kappa$ of $A$ along $\gamma= p \circ \vartheta$ is a Jacobi section along $\vartheta$ if there exist a dynamical variation $\zeta: \widetilde{I}\times I \rightarrow A$ of $\vartheta$ satisfying property (A) and a compatible transverse variation $\varkappa: \widetilde{I}\times I \rightarrow A$ of $\zeta$ such that $\varkappa(0,t)=\kappa(t)$.
\end{definition}

%%%%%%%%

\begin{remark}\label{rem:eq}
It follows from Theorem~\ref{theo-jacobi-transv} and \cite[Thm.~3.7]{carinena2015jacobi} that the Jacobi sections of
Definition~\ref{defnew-jacobi} are exactly the solutions of the Jacobi
equation $\nabla\nabla\kappa+\Phi(\kappa)=0$ of the semispray. In contrast
with the Jacobi equation of \cite{Avendano-Kama-2026,Boucetta2011}, which is formulated through
the curvature of a linear connection and therefore requires a spray, this
characterization holds for arbitrary semisprays.
\end{remark}
%\begin{theorem}
%Let $S\in \Gamma(\mathcal{T}A)$ be a SODE on $A$. The map% $ \Theta: \widetilde{I}\times I \rightarrow \mathcal{T}A$ 
%\begin{equation*}
%\Theta = S(\zeta(s,t))dt + B(s,t)ds.
%\end{equation*}
%with $B(s,t)=\left(\zeta(s,t),\varkappa(s,t),\frac{\partial \zeta}{\partial s}(s,t) \right) $ is a morphism of Lie algebroid over $\zeta: \widetilde{I}\times I \rightarrow A$ if and only if $\zeta$ is dynamical variation of $\vartheta(t):=\zeta(0,t)$ and $\varkappa$ is a compatible transverse variation of $\zeta$.
%\end{theorem}

\section{The anchor correspondence between Jacobi sections and Jacobi fields}

In this section, we clarify the relationship between Jacobi fields presented in Section 4, and the notion of  Jacobi sections revised in Section 6. The variational interpretation of Jacobi sections established in Theorem \ref{theo-jacobi-transv} will provide us with the appropriate framework to relate both notions. Indeed, while Jacobi fields arise from arbitrary dynamical variations of integral curves of a semispray, Jacobi sections correspond to those dynamical variations satisfying Property (A) and admitting compatible transverse variations. This additional geometric structure allows one to compare both notions through the anchor map. In this section, we establish this correspondence and characterize the circumstances under which a Jacobi field can be realized by a Jacobi section.

Our first result states the relation between Jacobi sections and Jacobi fields.

\begin{theorem}\label{sect-field-jac-theo}
Let $\mathcal{V}\in \mathfrak{X}(A)$ be a semispray on a Lie algebroid $A$, $\vartheta: I \rightarrow A$ an integral curve of $\mathcal{V}$ and $\gamma = p\circ \vartheta.$ If $\kappa(t)$ is a Jacobi section along $\vartheta$, then
\begin{equation}
J(t)= \rho(\kappa(t))
\end{equation} 
is a Jacobi field of $\mathcal{V}$ along $\gamma$.
\end{theorem}
\begin{proof}
 Suppose that $\kappa(t)$ is a Jacobi section along $\vartheta$. By Theorem \ref{theo-jacobi-transv}, there exists a dynamical variation $\zeta: \widetilde{I}\times I \rightarrow A$ of $\vartheta$ satisfying property (A) and a compatible transverse variation $\varkappa: \widetilde{I}\times I \rightarrow A$ of $\zeta$ such that $\varkappa(0,t)=\kappa(t)$. Then, $\alpha(s,t)= p \circ \zeta (s,t)$ is a variation of $\gamma$ by integral curves in the base. Thus, $\displaystyle J(t)=\frac{\partial \alpha}{\partial s}(0,t)$ is a Jacobi field along of $\gamma$. Since $\varkappa$ is a transverse variation, we have $\displaystyle  \frac{\partial \alpha}{\partial s}(s,t) = \rho(\varkappa(s,t))$. Therefore,
 $$ \rho(\kappa(t)) =  \rho(\varkappa(0,t))= \frac{\partial \alpha}{\partial s}(0,t) =J(t).$$
 \end{proof}

In summary, the previous theorem shows that a Jacobi section of a Lie algebroid is projected under the anchor map into a Jacobi field. However, the converse of the theorem holds under an additional hypothesis about the dynamical variations.

\begin{theorem}\label{sect-field-jac-theo2}
Let $\mathcal{V}\in \mathfrak{X}(A)$ be a semispray on a Lie algebroid $A$, $\vartheta: I \rightarrow A$ an integral curve of $\mathcal{V}$ and $\gamma = p\circ \vartheta.$ Suppose that $J(t)$ is a Jacobi field arising from a variation $\alpha: \tilde{I}\times I \rightarrow M$ by integral curves in the base whose image is completely contained in a leaf $\mathcal{L}  \subset M$ of the characteristic foliation, then there exists a Jacobi section $\kappa(t)$ along $\vartheta$ such that
\begin{equation}
J(t)= \rho(\kappa(t)).
\end{equation} 
\end{theorem}
\begin{proof}
Recall (Definition \ref{defjacfield} and Remark \ref{remark2}) that Jacobi fields are variational vector fields of dynamical variations in the base. By assumption, $J(t)=\frac{\partial \alpha}{\partial s}(0,t)$ for a variation $\alpha(s,t)=p\circ \zeta(s,t)$, where $\zeta: \tilde{I}\times I \rightarrow A$ is a dynamical variation of $\vartheta$. Since the image of $\alpha$ is completely contained in a leaf $\mathcal{L} $ of the characteristic foliation, $\zeta$ satisfies the property (A). By Lemma \ref{existencecontlift}, there is a curve $\varkappa_0: \tilde{I}\rightarrow A$ such that $\rho(\varkappa_0(s))= \frac{\partial\alpha}{\partial}(s,0).$ By Proposition \ref{varzero}, there exists a compatible transverse variation $\varkappa: \tilde{I}\times I \rightarrow A$ such that $\varkappa(s,0)= \varkappa_0(s)$. By Theorem  \ref{theo-jacobi-transv}, $\kappa(t)=\varkappa(0,t)$ is a Jacobi section along $\vartheta$. Finally, we have
 $$ \rho(\kappa(t)) =  \rho(\varkappa(0,t))= \frac{\partial \alpha}{\partial s}(0,t) =J(t).$$
\end{proof}

The theorems above completely describe the correspondence between Jacobi fields and Jacobi sections for arbitrary semisprays. We must emphasize that Jacobi sections always project through the anchor map to Jacobi fields. But not every Jacobi field comes, via the anchor, from some Jacobi sections. This phenomenon comes from the fact that the Jacobi fields arise from arbitrary dynamical variations of integral curves, whereas Jacobi sections arise from the additional geometric structure of dynamical variations satisfying Property (A) together with compatible transverse variations. Hence, the Jacobi field for semisprays on Lie algebroids presented in Section 4 is a geometric structure that considers a broader class of dynamical variations than the Jacobi sections.  
\vspace{12pt}

We conclude this section by comparing with the spray case treated in \cite{Avendano-Kama-2026}. In this paper, inspired by \cite{Boucetta2011}, Jacobi sections were defined as the sections along an integral curve satisfying a Jacobi equation written through the curvature of a linear connection. This is possible because a spray determines a unique torsionless linear connection whose geodesics are the integral curves of the spray \cite[Thm.~13]{Avendano-Kama-2026}. Without homogeneity, a semispray induces no linear connection, and this approach is no longer available.

For sprays,  Jacobi sections can also be obtained from both dynamical variations and compatible transverse variations
\cite[Prop.~17]{Avendano-Kama-2026}, \cite[Prop.~3.17]{Boucetta2011}. In Theorem~\ref{theo-jacobi-transv}, we show that this variational description does not depend on homogeneity: it holds for the Jacobi sections of \cite{carinena2015jacobi} of an arbitrary semispray. Homogeneity is only needed to write the Jacobi equation through a linear connection (compare Remark~\ref{rem:eq}).

%%%%%%%%

\section{An illustration of the correspondence between Jacobi sections and Jacobi fields}
In this section we present examples illustrating the variational framework developed in the previous sections. Our main goal is to show how compatible transverse variations arise naturally in concrete situations and how they provide the geometric mechanism relating dynamical variations, Jacobi sections and Jacobi fields. In particular, the examples illustrate the correspondence established in Section 7 and complement the abstract theory with explicit constructions.

Let $M=\mathbb{R}$ and the abelian Lie algebra $\mathfrak{g}=\mathbb{R}^2=\mathrm{span}\{ \mathbf{e}_1, \mathbf{e}_2\}$. Consider the infinitesimal action of $\mathfrak{g}$ on $\mathbb{R}$ defined by
\begin{equation*}
\Psi(\mathbf{e}_1 ) = x\frac{\partial}{\partial x }, \ \ \ \ \Psi(\mathbf{e}_2 ) =0.
\end{equation*}
In this example, we study the action Lie algebroid given by the vector bundle
$$ p: M\times  \mathbb{R}^2 \rightarrow M$$
with coordinates $(x,y,z)$ where the elements on the fiber througth $x$ is given by $a= y\mathbf{e}_1 + z\mathbf{e}_2$. The anchor map $\rho: M\times  \mathbb{R}^2 \rightarrow TM$ is
$$ \rho(x,y,z)= xy\frac{\partial}{\partial x}. $$
The bracket of sections can be defined in terms of the global basis sections $\sigma_1(x)=(x,1,0)$, $\sigma_2(x)=(x,0,1)$ by $ [ \sigma_1,\sigma_2 ]= 0, $ and extended to arbitrary sections by bilinearity, skew-symmetry and the Leibniz rule. The characteristic foliation has three components: $L_{-}=(-\infty,0)$, $L_0=\{ 0\},$ and $L_{+}=(0,\infty).$
Since $\rho(\sigma_1)=x\frac{\partial}{\partial x}$ and $\rho(\sigma_2)=0$, the semisprays have the following form
$$ \mathcal{V}(x,y,z) = xy\frac{\partial}{\partial x} +\nu_1(x,y,z)\frac{\partial}{\partial y}+\nu_2(x,y,z)\frac{\partial}{\partial z}.  $$
For simplicity, we consider the semispray corresponding to  $\nu_1=\nu_2 \equiv 0$. The corresponding integral curves satisfy the following system
$$ \dot{x}=xy,\ \ \ \ \dot{y}=0,\ \ \ \ \dot{z}=0. $$
Hence, the solution through the point $(x_0,y_0,z_0)$ at $t=0$ is the curve $\vartheta: \mathbb{R} \rightarrow M\times  \mathbb{R}^2 $
$$ \vartheta(t)=(x_0e^{y_0t},y_0,z_0) $$
with $\gamma(t)=x_0e^{y_0t}$. If $x_0>0$, then the dynamical variation $\zeta: \mathbb{R}\times \mathbb{R} \rightarrow  M\times  \mathbb{R}^2 $ of $\vartheta$
$$\zeta(s,t)= (x_0e^{(y_0+rs)t+us},y_0+rs,z_0+vs), $$
where $r,u,v$ are real constants, satisfies the property (A), because the image $\alpha(s,t)= x_0e^{(y_0+rs)t+us}$ is contained in $L_{+}$. By Theorem \ref{sect-field-jac-theo2}, the Jacobi field along $\gamma$
$$J(t)=\frac{\partial \alpha}{\partial s}(0,t)=\left(x_0(rt+u)e^{y_0t} \right)\frac{\partial }{\partial x}, $$
is in the image of the anchor map. On the other hand, the map $\varkappa: \mathbb{R}\times \mathbb{R} \rightarrow  M\times  \mathbb{R}^2 $ defined by
$$\varkappa(s,t)=(x_0e^{(y_0+rs)t+us}, rt+u, vt )  $$
is a transverse variation of $\zeta$. Using equation \eqref{lineardiffeq}, we can see that $\varkappa$ is also compatible. So, $\kappa(t)=(x_0e^{y_0t}, u+rt,vt)$ is a Jacobi section along $\vartheta$ and
$$ \rho(\kappa(t))= J(t). $$
Therefore, this example provides an explicit illustration of the correspondence established in Theorem  \ref{sect-field-jac-theo2}: the Jacobi field $J(t)$ arises as the anchor image of the Jacobi section $\kappa(t)$, obtained from a compatible transverse variation.

To emphasize that Property (A) is an essential hypothesis in Theorem \ref{sect-field-jac-theo2}, we now consider a dynamical variation that does not satisfy this property. The map
$$\widetilde{\zeta}(s,t)= (se^{(y_0+rs)t)},y_0+rs,z_0+us ) $$ 
is a dynamical variation of the integral curve $\widetilde{\vartheta}(t)=(0,y_0,z_0)$. Observe that the image of $\widetilde{\alpha}(s,t)=se^{(y_0+rs)t)}$ is not contained in any single characteristic leaf, since it intersects $L_{-}$, $L_0,$ and $L_{+}.$ The Jacobi field associated with this variation is
$$\widetilde{J}(t)=e^{y_0t} \frac{\partial }{\partial x}$$
Since $\widetilde{\gamma}(t)=p\circ \widetilde{\vartheta}(t)\equiv 0$, any section along $\widetilde{\gamma}$ has the following form $\varsigma(t)=(0, \varsigma_1(t),\varsigma_2(t))$. Hence, there is no Jacobi section $\kappa(t)$ along $\widetilde{\vartheta}(t)$ such that $\rho(\kappa(t))=\widetilde{J}(t)$. Although $\widetilde{\alpha}$ generates a Jacobi field, its image is not contained in a single characteristic leaf, and therefore no compatible transverse variation can exist, and hence no Jacobi section can be associated with it.

\paragraph{Acknowledgements.} 
SECIHTI supported this research under the Grant CF-2023-G-279.
 
%martinez2002lie  carinena2007reduction

\end{document}